# Integrable Hamiltonian systems on Lie groups: Kowalewski type

By V. Jurdjevic

### Introduction

The contributions of Sophya Kowalewski to the integrability theory of the equations for the heavy top extend to a larger class of Hamiltonian systems on Lie groups; this paper explains these extensions, and along the way reveals further geometric significance of her work in the theory of elliptic curves. Specifically, in this paper we shall be concerned with the solutions of the following differential system in six variables $h_1, h_2, h_3, H_1, H_2, H_3$

$$\frac{dH_1}{dt} = H_2 H_3 \left( \frac{1}{c_3} - \frac{1}{c_2} \right) + h_2 a_3 - h_3 a_2 \ ,$$

$$\frac{dH_2}{dt} = H_1 H_3 \left( \frac{1}{c_1} - \frac{1}{c_3} \right) + h_3 a_1 - h_1 a_3 \ ,$$

$$\frac{dH_3}{dt} = H_1 H_2 \left( \frac{1}{c_2} - \frac{1}{c_1} \right) + h_1 a_2 - h_2 a_1 \ ,$$

$$\frac{dh_1}{dt} = \frac{h_2 H_3}{c_3} - \frac{h_3 H_2}{c_2} + k(H_2 a_3 - H_3 a_2) \ ,$$

$$\frac{dh_2}{dt} = \frac{h_3 H_1}{c_1} - \frac{h_1 H_3}{c_3} + k(H_3 a_1 - H_1 a_3) \ ,$$

$$\frac{dh_3}{dt} = \frac{h_1 H_2}{c_2} - \frac{h_2 H_1}{c_1} + k(H_1 a_2 - H_2 a_1) \ ,$$

in which $a_1, a_2, a_3, c_1, c_2, c_3$ and $k$ are constants. The preceding system of equations can also be written more compactly

(i) $$\frac{d\widehat{H}}{dt} \ = \ \widehat{H} \times \widehat{\Omega} + \hat{h} \times \hat{a}, \quad \frac{d\hat{h}}{dt} \ = \ \hat{h} \times \widehat{\Omega} + k(\widehat{H} \times \hat{a})$$

with $\times$ denoting the vector product in $\mathbb{R}^3$ and with

$$\widehat{H} \ = \ \begin{pmatrix} H_1 \\ H_2 \\ H_3 \end{pmatrix}, \quad \widehat{\Omega} \ = \ \begin{pmatrix} \frac{H_1}{c_1} \\ \frac{H_2}{c_2} \\ \frac{H_3}{c_3} \end{pmatrix}, \quad \hat{h} \ = \ \begin{pmatrix} h_1 \\ h_2 \\ h_3 \end{pmatrix} \quad \text{and} \quad \hat{a} \ = \ \begin{pmatrix} a_1 \\ a_2 \\ a_3 \end{pmatrix}.$$

When $k = 0$ the preceding equations formally coincide with the equations of the motions of a rigid body around its fixed point in the presence of the gravitational force, known as the heavy top in the literature on mechanics. In this context, the constants $c_1, c_2, c_3$ correspond to the principal moments of inertia of the body, while $a_1, a_2, a_3$ correspond to the coordinates of the center of mass of the body relative to an orthonormal frame fixed on the body, known as the moving frame. The vector $\Omega$ corresponds to the angular velocity of the body measured relative to the moving frame. That is, if $R(t)$ denotes the orthogonal matrix describing the coordinates of the moving frame with respect to a fixed orthonormal frame, then

$$\frac{dR(t)}{dt} \; = \; R(t) \begin{pmatrix} 0 & -\Omega_3(t) & \Omega_2(t) \\ \Omega_3(t) & 0 & -\Omega_1(t) \\ -\Omega_2(t) & \Omega_1(t) & 0 \end{pmatrix} \; .$$

The vector $\widehat{H}$ corresponds to the angular momentum of the body, related to the angular velocity by the classic formulas $\frac{1}{c_i}H_i = \Omega_i$, $i = 1, 3, 3$. Finally, the vector $\hat{h}(t)$ corresponds to the movements of the vertical unit vector observed from the moving body and is given by, $\hat{h}(t) = R^{-1}(t) \begin{pmatrix} 0 \\ 0 \\ 1 \end{pmatrix}$. Therefore, solutions of equations (i) corresponding to $k = 0$, and further restricted to $h_1^2 + h_2^2 + h_3^2 = 1$ coincide with all possible movements of the heavy top.

Rather than studying the foregoing differential system in $\mathbb{R}^6$, as is commonly done in the literature of the heavy top, we shall consider it instead as a Hamiltonian system on the group of motions $E_3$ of a Euclidean space $E^3$ corresponding to the Hamiltonian function

(ii)         $$H \; = \; \frac{1}{2}\left( \frac{H_1^2}{c_1} + \frac{H_2^2}{c_2} + \frac{H_3^2}{c_3} \right) + a_1 h_1 + a_2 h_2 + a_3 h_3.$$

This Hamiltonian system has its origins in a famous paper of Kirchhoff of 1859 concerning the equilibrium configurations of an elastic rod, in which he likened the basic equations of the rod to the equations of the heavy top. His observation has since been known as the kinetic analogue of the elastic rod. According to Kirchhoff an elastic rod is modeled by a curve $\gamma(t)$ in a Euclidean space $E^3$ together with an orthonormal frame defined along $\gamma(t)$ and adapted to the curve in a prescribed manner. The usual assumptions are that the rod is inextensible, and therefore $\|\frac{d\gamma}{dt}\| = 1$, and that the first leg of the frame coincides with the tangent vector $\frac{d\gamma}{dt}$. In this context, $\gamma(t)$ corresponds to the central line of the rod, and the frame along $\gamma$ measures the amount of bending and twisting of the rod relative to a standard reference frame defined by the unstressed state of the rod.

Denote by $R(t)$ the relation of the frame along $\gamma$ to the reference frame; then $R(t)$ is a curve in $\mathrm{SO}_3(\mathbb{R})$ and therefore

$$\frac{dR(t)}{dt} = R(t) \begin{pmatrix} 0 & -u_3(t) & u_2(t) \\ u_3(t) & 0 & -u_1(t) \\ -u_2(t) & u_1(t) & 0 \end{pmatrix}$$

for functions $u_1(t), u_2(t), u_3(t)$. In the literature on elasticity these functions are called strains. Kirchhoff's model for the equilibrium configurations of the rod subject to the prescribed boundary conditions, consisting of the terminal positions of the rod and its initial and final frame, postulates that the equilibrium configurations minimize the total elastic energy of the rod $\frac{1}{2} \int_0^T (c_1 u_1^2(t) + c_2 u_2^2(t) + c_3 u_3^2(t)) dt$ with $c_1, c_2, c_3$ constants, determined by the physical characteristics of the rod, with $T$ equal to the length of the rod.

From the geometric point of view each configuration of the rod is a curve in the frame bundle of $E^3$ given by the following differential system

(iii) $\qquad \dfrac{d\gamma}{dt} = R(t) \begin{pmatrix} a_1 \\ a_2 \\ a_3 \end{pmatrix}, \qquad \dfrac{dR}{dt} = R(t) \begin{pmatrix} 0 & -u_3 & u_2 \\ u_3 & 0 & -u_1 \\ -u_2 & u_1 & 0 \end{pmatrix}$

with constants $a_1, a_2, a_3$ describing the relation of the tangent vector $\frac{d\gamma}{dt}$ to the frame along $\gamma$. The preceding differential system has a natural interpretation as a differential system in the group of motions $E_3 = E^3 \rtimes \mathrm{SO}_3(\mathbb{R})$. The Hamiltonian $H$ given above appears as a necessary condition of optimality for the variational problem of Kirchhoff.

In contrast to the traditional view of applied mathematics influenced by Kirchhoff, in which the elastic problem is likened to the heavy top, we shall show that the analogy goes the other way; the heavy top is like the elastic problem and much of the understanding of the integrability of its basic equations is gained through this analogy. To begin with, the elastic problem, dependent only on the Riemannian structure of the ambient space extends to other Riemannian spaces. In particular, for spaces of constant curvature, the frame bundle is identified with the isometry group, and the parameter $k$ that appears in the above differential system coincides with their curvature. In this paper we shall concentrate on $k = 0$, $k = \pm 1$. The case $k = 1$, called the elliptic case, corresponds to the sphere $S^3 = \mathrm{SO}_4(\mathbb{R})/\mathrm{SO}_3(\mathbb{R})$, while $k = -1$, called the hyperbolic case, corresponds to the hyperboloid $\mathbb{H}^3 = \mathrm{SO}(1,3)/\mathrm{SO}_3(\mathbb{R})$. As will be shown subsequently, differential systems described by (i) correspond to the projections of Hamiltonian differential equations on the Lie algebra of $G$ generated by the Hamiltonian $H$ in (ii), with $G$ any of the groups $E_3, \mathrm{SO}_4(\mathbb{R})$ and $\mathrm{SO}(1,3)$ as the isometry groups of the above symmetric spaces.

It may be relevant to observe that equations (iii) reduce to Serret-Frenet equations for a curve $\gamma$ when $u_2 = 0$. Then $u_1(t)$ is the torsion of $\tau(t)$ of $\gamma$ while $u_3$ is its curvature $\kappa(t)$. Hence the elastic energy of $\gamma$ becomes a functional of its geometric invariants. In particular, the variational problem attached to $\int_0^T (\kappa^2(t) + \tau^2(t)) dt$ was considered by P. Griffith a natural candidate for the elastic energy of a curve. Equations (i) then correspond to this variational problem when $a_2 = a_3 = 1$, $c_1 = c_3 = 1$, and $c_2 = \infty$.

With these physical and geometric origins in mind we shall refer to this class of Hamiltonian systems as elastic, and refer to the projections of the integral curves of the corresponding Hamiltonian vector field on the underlying symmetric space as elastic curves. Returning now to our earlier claim that much of the geometry of the heavy top is clarified through the elastic problem, we note that, in contrast to the heavy top, the elastic problem is a left-invariant variational problem on $G$, and consequently always has five independent integrals of motion.

These integrals of motion are $H$ itself, two Casimir integrals $\|\hat{h}\|^2 + k\|\widehat{H}\|^2$, $\hat{h} \cdot \widehat{H} = h_1 H_1 + h_2 H_2 + H_3 h_3$, and two additional integrals due to left-invariant symmetry determined by the rank of the Lie algebra of $G$. This observation alone clarifies the integrability theory of the heavy top as it demonstrates that the existence of a fourth integral for differential system (i) is sufficient for its complete integrability.

It turns out that completely integrable cases for the elastic problem occur under the same conditions as in the case of the heavy top. In particular, we have the following cases:

(1) $a = 0$. Then, both $\|\hat{h}\|$ and $\|\widehat{H}\|$ are integrals of motion. This case corresponds to *Euler's top*. The elastic curves are the projections of the extremal curves in the intersection of energy ellipsoid $H = \frac{1}{2}\left(\frac{H_1^2}{c_1} + \frac{H_2^2}{c_2} + \frac{H_3^2}{c_3}\right)$ with the momentum sphere $M = H_1^2 + H_2^2 + H_3^2$.

(2) $c_2 = c_3$ and $a_2 = a_3 = 0$. In this case $H_1$ is also an integral of motion. This case corresponds to *Lagrange's top*. Its equations are treated in complete detail in [8].

(3) $c_1 = c_2 = c_3$. Then $H_1 a_1 + H_2 a_2 + H_3 a_3$ is an integral of motion. This integral is also well-known in the literature of the heavy top. The corresponding equations are integrated by means of elliptic functions similar to the case of Lagrange, which partly accounts for its undistinguished place in the hierarchy of integrable tops.

The remaining, and the most fascinating integrable case was discovered by S. Kowalewski in her famous paper of 1889 under the conditions that $c_1 = c_2 = 2c_3$ and $a_3 = 0$. It turns out that the extra integral of motion exists under the same conditions for the elastic problem, and is equal to

$$|z^2 - a(w - ka)|^2$$

with $z = \frac{1}{2}(H_1 + iH_2)$, $w = h_1 + ih_2$ and $a = a_1 + ia_2$. This integral formally coincides with that found by Kowalewski only for $k = 0$.

The present paper is essentially devoted to this case. We shall show that Kowalewski's method of integration extends to the elastic problem with only minor modifications and leads to hyperelliptic differential equations on Abelian varieties on the Lie algebra of G. Faced with the "mysterious change of variables" in Kowalewski's paper, whose mathematical nature was never properly explained in the literature of the heavy top, we discovered simple and direct proofs of the main steps that not only clarify Kowalewski's method but also identify Hamiltonian systems as an important ingredient of the theory of elliptic functions.

As a byproduct this paper offers an elementary proof of Euler's results of 1765 concerning the solutions of

$$\frac{dx}{\sqrt{P(x)}} \pm \frac{dy}{\sqrt{P(y)}} = 0$$

with $P$ an arbitrary fourth degree polynomial with complex coefficients. Combined further with A. Weil's interpretations of Euler's results in terms of addition formulas for curves $u^2 = P(x)$, these results form a theoretic base required for the integration of the extremal equations.

This seemingly unexpected connection between Kowalewski, Euler and Weil is easily explained as follows:

The elastic problem generates a polynomial $P(x)$ of degree four and two forms $R(x, y)$ and $\widehat{R}(x, y)$ each of degree four satisfying the following relations

(iv) $\quad R(x, x) = P(x)$, and $\quad R^2(x, y) + (x - y)^2 \widehat{R}(x, y) = P(x)P(y)$.

We begin our investigations with these relations associated to an arbitrary polynomial $P(x) = A + 4Bx + 6Cx^2 + 4Dx^3 + Ex^4$. In particular, we explicitly calculate the coefficients of $\widehat{R}$ corresponding to $R(x, y) = A + 2B(x + y) + 3C(x^2 + y^2) + 2Dxy(x + y) + Ex^2y^2$. Having obtained the expression for $\widehat{R}$, we have $\widehat{R}_\theta(x, y) = -(x - y)^2\theta^2 + 2R(x, y)\theta + \widehat{R}(x, y)$ is the form in (iv) that corresponds to the most general form $R_\theta(x, y) = R(x, y) - \theta(x - y)^2$, that satisfies $R(x, x) = P(x)$.

We then show that $\widehat{R}_\theta(x, y) = 0$ contains all solutions of $\frac{dx}{\sqrt{P(x)}} \pm \frac{dy}{\sqrt{P(y)}} = 0$ as $\theta$ varies over all complex numbers. This demonstration recovers the result of Euler in 1765, and also identifies the parametrizing variable $\theta$ with the points on the canonical cubic elliptic curve

$$\Gamma = \{(\xi, \eta): \eta^2 = 4\xi^3 - g_2\xi - g_3\}$$

with

$$g_2 = AE - 4BD + 3C^2 \quad \text{and} \quad g_3 = ACE + 2BCD - AD^2 - B^2E - C^3$$

via the relation $\theta = 2(\xi + C)$. The constants $g_2$ and $g_3$ are known as the covariant invariants of the elliptic curve $\mathcal{C} = \{(x, y) : u^2 = P(x)\}$. André Weil points out in [13] that the results of Euler have algebraic interpretations that may be used to define an algebraic group structure on $\Gamma \cup C$. The "mysterious" change of variables in the paper of Kowalewski is nothing more than the transformation from $\mathcal{C} \times \mathcal{C}$ into $\Gamma \times \Gamma$ given by $(N, M) \in \mathcal{C} \times \mathcal{C}$ to $(O_1, O_2)$ in $\Gamma \times \Gamma$ with $O_1 = N - M$ and $O_2 = N + M$.

The actual formulas that appear in the paper of Kowalewski

$$\frac{d\xi_1}{\eta(\xi_1)} = -\frac{dx}{\sqrt{P(x)}} + \frac{dy}{\sqrt{P(y)}} \quad \text{and} \quad \frac{d\xi_2}{\eta(\xi_2)} = \frac{dx}{\sqrt{P(x)}} + \frac{dy}{\sqrt{P(y)}}$$

are the infinitesimal versions of Weil's addition formulas.

Oddly enough, Kowalewski omits any explanation concerning the origins and the use of the above formulas, although it seems very likely that the connections with the work of Euler were known to her at that time (possibly through her association with Weierstrass).

The organization of this paper is as follows: Section I contains a self-contained treatment of Hamiltonian systems on Lie groups. This material provides a theoretic base for differential equations (i) and their conservation laws. In contrast to the traditional treatment of this subject matter geared to the applications in mechanics, the present treatment emphasizes the geometric nature of the subject seen through the left-invariant realization of the symplectic form on $T^*G$, the latter considered as $G \times \mathfrak{g}^*$ via the left-translations.

Section II contains the reductions in differential equations through the conservation laws (integrals of motion) leading to the fundamental relations that appear in the paper of Kowalewski. Section III contains the proof of Euler's result along with its algebraic interpretations by A. Weil. Section IV explains the procedure for integrating the differential equations by quadrature leading up to the famous hyperelliptic curve of Kowalewski.

Section V is devoted to complex extensions of differential system (i) in which the time variable is also considered complex. Motivated by the brilliant observation of Kowalewski that integrable cases of the heavy top are integrated by means of elliptic and hyperelliptic integrals and that, therefore, the solutions are meromorphic functions of complex time, we investigate the cases of elastic equations that admit purely meromorphic solutions on at least an open subset of $\mathbb{C}^6$ under the assumption that $c_1 = c_2$, while the remaining coefficient $c_3$ is arbitrary. We confirm Kowalewski's claim in this more general setting that the only cases that admit such meromorphic solutions are the ones already described in our introduction. In doing so we are unfortunately obliged to make an additional assumption (that is likely inessential) concerning the order of poles in solutions. This assumption is necessitated by a gap in Kowalewski's original paper, first noticed by A.A. Markov, that apparently still remains open in the literature on the heavy top. We conclude the paper with an integrable (in the sense of Liouville) elastic case that falls outside of the meromorphically integrable class suggesting further limitations of Kowalewski's methods in the classification of completely integrable elastic systems.

## 1. Hamiltonian systems on Lie groups

We shall use $\mathfrak{g}$ to denote the Lie algebra of a Lie group $G$, while $\mathfrak{g}^*$ will denote the dual of $\mathfrak{g}$. The cotangent bundle $T^*G$ will be identified with $G \times \mathfrak{g}^*$ via the left-translations: an element $(g, p)$ in $G \times \mathfrak{g}^*$ is identified with $\xi \in T_g^*G$ by $p = dL_g^* \xi$ with $dL_g^*$ denoting the pull-back of the left-translation $L_g(x) = gx$.

The tangent bundle of $T^*G$ is identified with $TG \times \mathfrak{g}^* \times \mathfrak{g}^*$, the latter further identified with $G \times \mathfrak{g}^* \times \mathfrak{g} \times \mathfrak{g}^*$. Relative to this decomposition, vector fields on $T^*G$ will be denoted by $(X(g, p), Y^*(g, p))$ with $(g, p)$ denoting the base point in $T^*G$ and $X$ and $Y^*$ denoting their values in $\mathfrak{g}$ and $\mathfrak{g}^*$ respectively. Then the canonical symplectic form $\omega$ on $T^*G$ in the aforementioned trivialization of $T^*G$ is given by:

$$(1) \qquad \omega_{(g,p)}\left((X_1, Y_1^*), (X_2, Y_2^*)\right) \;=\; Y_2^*(X_1) - Y_1^*(X_2) - p[X_1, X_2] \;.$$

The correct signs in this expression depend on the particular choice of the Lie bracket. For the above choice of signs, $[X, Y](f) = Y(Xf) - X(Yf)$ for any function $f$.

The symplectic form sets up a correspondence between functions $H$ on $T^*G$ and vector fields $\vec{H}$ given by

$$(2) \qquad \omega_{(g,p)}(\vec{H}(g, p), V) \;=\; dH_{(g,p)}(V)$$

for all tangent vectors $V$ at $(g, p)$. It is customary to call $H$ a Hamiltonian function, or simply a Hamiltonian, and $\vec{H}$ the Hamiltonian vector field of $H$. A Hamiltonian $H$ is called left-invariant if it is invariant under left-translations, which is equivalent to saying that $H$ is a function on $\mathfrak{g}^*$; that is, $H$ is constant over the fiber above each point $p$ in $\mathfrak{g}^*$.

For each left-invariant Hamiltonian $H$, $dH$, being a linear function over $\mathfrak{g}^*$, is an element of $\mathfrak{g}$ at each point $p$ in $\mathfrak{g}^*$. Then it follows from (1) and (2) that the Hamiltonian vector field $\vec{H}$ of a left-invariant Hamiltonian is given by $(X(p), Y^*(p))$ in $\mathfrak{g} \times \mathfrak{g}^*$ with

$$(3) \qquad X(p) \;=\; dH_p \quad \text{and} \quad Y^*(p) \;=\; -\mathrm{ad}^*(dH_p)(p).$$

In this expression $\mathrm{ad}^*X$ denotes the dual mapping of $\mathrm{ad}X : \mathfrak{g} \to \mathfrak{g}$ given by $\mathrm{ad}X(Y) = [X, Y]$.

It immediately follows from (3) that the integral curves $(g(t), p(t))$ of $\vec{H}$ satisfy

$$(4) \qquad dL_{g^{-1}(t)} \frac{dg}{dt} \; = \; dH_{p(t)} \quad \text{and} \quad \frac{dp}{dt} \; = \; -\text{ad}^*(dH_{p(t)})(p(t)) \; ,$$

and consequently

$$(5) \qquad\qquad p(t) \; = \; \text{Ad}^*_{g(t)}(p(0))$$

with $\text{Ad}^*$ equal to the co-adjoint action of $G$ on $\mathfrak{g}^*$. Thus, the projections of integral curves of left-invariant Hamiltonian vector fields evolve on the co-adjoint orbits of $G$.

When the group $G$ is semisimple the Killing form is nondegenerate and can be used to identify elements in $\mathfrak{g}^*$ with elements in $\mathfrak{g}$. This correspondence identifies each curve $p(t)$ in $\mathfrak{g}^*$ with a curve $U(t)$ in $\mathfrak{g}$. For integral curves of a left-invariant Hamiltonian $H$, the equation $\frac{dp}{dt} \; = \; -\text{ad}^*(dH_p)(p(t))$ corresponds to

$$(6) \qquad\qquad \frac{dU(t)}{dt} \; = \; [dH_{p(t)}, U(t)] \; .$$

The expression (6) is often called the Lax-pair form in the literature on the Hamiltonian systems.

We shall use $\{F, H\}$ to denote the Poisson bracket of functions $F$ and $H$. Recall that $\{F, H\}(g, p) = \omega_{(g,p)}\{\vec{F}(g), \vec{H}(p)\}$. It follows immediately from (1) that for left-invariant Hamiltonians $F$ and $H$, their Poisson bracket is given by $\{F, H\}(p) = p([dF_p, dH_p])$, for all $p$ in $\mathfrak{g}^*$.

A function $F$ on $T^*G$ is called an integral of motion for $H$ if $F$ is constant along each integral curve of $\vec{H}$, or equivalently if $\{F, H\} = 0$. A given Hamiltonian is said to be completely integrable if there exist $n - 1$ independent integrals of motion $F_1, \dots, F_{n-1}$ that together with $F_n = H$ satisfy $\{F_i, F_j\} = 0$ for all $i, j$. The independence of $F_1, \dots, F_n$ is taken in the sense that the differentials $dF_1, \dots, dF_n$ are independent at all points on $T^*G$.

Any vector field $X$ on $G$ lifts to a function $F_X$ on $T^*G$ defined by $F_X(\xi) = \xi(X(g))$ for any $\xi \in T^*_g G$. In the left-invariant representation $G \times \mathfrak{g}^*$, left invariant vector fields lift to linear functions on $\mathfrak{g}^*$, while right-invariant vector fields lift to $F_X(g, p) = p(dL_{g^{-1}} \circ dR_g X_e)$ with $X_e$ denoting the value of $X$ at the group identity $e$ of $G$. The preceding expression for $F_X$ can also be written as $F_X(g, p) \; = \; \text{Ad}^*_{g^{-1}}(p)(X_e)$. Therefore, along each integral curve $(g(t), p(t))$ of a left-invariant Hamiltonian $\vec{H}$,

$$\begin{aligned} F_X(g(t), p(t)) \; &= \; \text{Ad}^*_{g^{-1}(t)} p(t) X_e \\ &= \; \text{Ad}^*_{g^{-1}(t)} \circ \text{Ad}^*_{g(t)}(p(0)) X_e \; = \; (p(0)) X_e \end{aligned}$$

and consequently $F_X$ is an integral of motion for $H$.

The maximum number of right-invariant vector fields that pairwise commute with each other is equal to the rank of $\mathfrak{g}$. Consequently, a left-invariant Hamiltonian $H$ always has $r$-independent integrals of motion Poisson commuting with each other, and of course commuting with $H$, with $r$ equal to the rank of $\mathfrak{g}$.

In addition to these integrals of motion, there may be functions on $\mathfrak{g}^*$ that are constant on co-adjoint orbits of $G$. Such functions are called Casimir functions, and they are integrals of motion for any left-invariant Hamiltonian $H$. On semisimple Lie groups Casimir functions always exist as can be seen from the Lax-pair representation (6). They are the coefficients of the characteristic polynomial of $U(t)$ (realized as a curve on the space of matrices via the adjoint representation).

With these concepts and this notation at our disposal we shall take $\mathfrak{g}$ to be any six dimensional Lie algebra with a basis $B_1, B_2, B_3, A_1, A_2, A_3$ that satisfies the following Lie bracket table:

| [ , ] | $A_1$ | $A_2$ | $A_3$ | $B_1$ | $B_2$ | $B_3$ |
|---|---|---|---|---|---|---|
| $A_1$ | $0$ | $-A_3$ | $A_2$ | $0$ | $-B_3$ | $B_2$ |
| $A_2$ | $A_3$ | $0$ | $-A_1$ | $B_3$ | $0$ | $-B_1$ |
| $A_3$ | $-A_2$ | $A_1$ | $0$ | $-B_2$ | $B_1$ | $0$ |
| $B_1$ | $0$ | $-B_3$ | $B_2$ | $0$ | $-kA_3$ | $kA_2$ |
| $B_2$ | $B_3$ | $0$ | $-B_1$ | $kA_3$ | $0$ | $-kA_1$ |
| $B_3$ | $-B_2$ | $B_1$ | $0$ | $-kA_2$ | $kA_1$ | $0$ |

with $\quad k = \left\{ \begin{array}{c} 0 \\ 1 \\ -1 \end{array} \right.$.

Table 1

The reader may easily verify that the following six dimensional matrices

$$B_1 = \begin{pmatrix} 0 & -k & 0 & 0 \\ 1 & 0 & 0 & 0 \\ 0 & 0 & 0 & 0 \\ 0 & 0 & 0 & 0 \end{pmatrix}, \ A_1 = \begin{pmatrix} 0 & 0 & 0 & 0 \\ 0 & 0 & 0 & 0 \\ 0 & 0 & 0 & -1 \\ 0 & 0 & 1 & 0 \end{pmatrix},$$

$$B_2 = \begin{pmatrix} 0 & 0 & -k & 0 \\ 0 & 0 & 0 & 0 \\ 1 & 0 & 0 & 0 \\ 0 & 0 & 0 & 0 \end{pmatrix}, \ A_2 = \begin{pmatrix} 0 & 0 & 0 & 0 \\ 0 & 0 & 0 & 1 \\ 0 & 0 & 0 & 0 \\ 0 & -1 & 0 & 0 \end{pmatrix},$$

$$B_3 = \begin{pmatrix} 0 & 0 & 0 & -k \\ 0 & 0 & 0 & 0 \\ 0 & 0 & 0 & 0 \\ 1 & 0 & 0 & 0 \end{pmatrix}, A_3 = \begin{pmatrix} 0 & 0 & 0 & 0 \\ 0 & 0 & -1 & 0 \\ 0 & 1 & 0 & 0 \\ 0 & 0 & 0 & 0 \end{pmatrix}$$

satisfy the above Lie bracket table under the matrix commutator bracket $[M, N] = NM - MN$. For $k = 0$, $\mathfrak{g}$ is the semi-direct product $\mathbb{R}^3 \ltimes \mathrm{so}_3(R)$, for $k = 1$, $\mathfrak{g}$ is $\mathrm{so}_4(R)$, and for $k = -1$, $\mathfrak{g} = \mathrm{so}(1,3)$.

Throughout this paper we shall use $h_i$ and $H_i$, to denote the linear functions on $\mathfrak{g}^*$ given by $h_i(p) = p(B_i)$, and $H_i(p) = p(A_i)$, $i = 1, 2, 3$. These functions are Hamiltonian lifts of left-invariant vector fields induced by the above basis in $\mathfrak{g}$. Finally, as stated earlier, we shall consider a fixed Hamiltonian function $H$ on $\mathfrak{g}^*$ given by

$$H = \frac{1}{2}\left(\frac{H_1^2}{c_1} + \frac{H_2^2}{c_2} + \frac{H_3^2}{c_3}\right) + a_1 h_1 + a_2 h_2 + a_3 h_3$$

for some constants $c_1, c_2, c_3$ and $a_1, a_2, a_3$.

We shall refer to the integral curves of $\vec{H}$ as the *extremal curves*. For each extremal curve $(g(t), p(t))$, $x(t) = g(t)e_1$ will be called an *elastic curve*. Elastic curves are the projections of the extremal curves on the underlying symmetric space $G/K$ with $K$ denoting the group that stabilizes $e_1$ in $\mathbb{R}^4$ (written as the column vector, with the action coinciding with the matrix multiplication). It can be easily verified that $K \simeq \mathrm{SO}_3(R)$ and that $G/K$ is equal to $\mathbb{R}^3, S^3$ or $\mathbb{H}^3$ depending whether $k = 0, 1$ or $-1$. The remaining columns of $g$ give the coordinates of the moving frame $v_1(t), v_2(t), v_3(t)$ defined along $x(t)$, and adapted to the curve $x(t)$ so that $\frac{dx(t)}{dt} = a_1 v_1(t) + a_2 v_2(t) + a_3 v_3(t)$.

*The semi-simple case.* For $k \neq 0$, the Killing form $T$ is nondegenerate and invariant in the sense that $T([A, B], C) = T(A, [B, C])$. We shall take $T(A, B) = \frac{1}{2}\mathrm{trace}\,(AB)$. It follows that $T(A, B) = -(\sum_{i=1}^{3} a_i \bar{a}_i + k b_i \bar{b}_i)$, with $A = \sum_{i=1}^{3} a_i A_i + b_i B_i$ and $B = \sum_{i=1}^{3} \bar{a}_i A_i + \bar{b}_i B_i$. Upon identifying $p$ in $\mathfrak{g}^*$ with $U$ in $\mathfrak{g}$ via the trace form, we get that

$$(7) \qquad U = \begin{pmatrix} 0 & h_1 & h_2 & h_3 \\ -kh_1 & 0 & H_3 & -H_2 \\ -kh_2 & -H_3 & 0 & H_1 \\ -kh_3 & H_2 & -H_1 & 0 \end{pmatrix}.$$

Then

$$dH_p = \begin{pmatrix} 0 & -ka_1 & -ka_2 & -ka_3 \\ a_1 & 0 & -\frac{1}{c_3}H_3(p) & \frac{1}{c_2}H_2(p) \\ a_2 & \frac{1}{c_3}H_3(p) & 0 & -\frac{1}{c_1}H_1(p) \\ a_3 & -\frac{1}{c_2}H_2(p) & \frac{1}{c_1}H_1(p) & 0 \end{pmatrix}$$

and equation (6) yields the following differential system:

$$\text{(8)} \qquad \frac{dh_1}{dt} = \frac{h_2 H_3}{c_3} - \frac{h_3 H_2}{c_2} + k(H_2 a_3 - H_3 a_2) \; ,$$

$$\frac{dh_2}{dt} = \frac{h_3 H_1}{c_1} - \frac{h_1 H_3}{c_3} + k(H_3 a_1 - H_1 a_3) \; ,$$

$$\frac{dh_3}{dt} = \frac{h_1 H_2}{c_2} - \frac{h_2 H_1}{c_1} + k(H_1 a_2 - H_2 a_1) \; ,$$

$$\frac{dH_1}{dt} = \frac{H_2 H_3}{c_3} - \frac{H_2 H_3}{c_2} + (h_2 a_3 - h_3 a_2) \; ,$$

$$\frac{dH_2}{dt} = \frac{H_1 H_3}{c_1} - \frac{H_1 H_3}{c_3} + (h_3 a_1 - h_1 a_3) \; ,$$

$$\frac{dH_3}{dt} = \frac{H_1 H_2}{c_2} - \frac{H_1 H_2}{c_1} + (h_1 a_2 - h_2 a_1) \; .$$

*Remark.* The foregoing differential equation can also be obtained by using the Poisson bracket through the formulas

$$\frac{dh_i}{dt} \; = \; \{h_i, H\}, \quad \text{and} \quad \frac{dH_i}{dt} \; = \; \{H_i, H\}, \quad i = 1, 2, 3.$$

Apart from the vector product representation given by (i) of the introduction, differential system (8) has several other representations. The most immediate, that will be useful for Section V, is the representation in $so_3(R) \times so_3(R)$ obtained by identifying vectors $\widehat{A} = \begin{pmatrix} \alpha_1 \\ \alpha_2 \\ \alpha_3 \end{pmatrix}$ in $\mathbb{R}^3$ with antisymmetric matrices $A = \begin{pmatrix} 0 & -\alpha_3 & \alpha_2 \\ \alpha_3 & 0 & -\alpha_1 \\ -\alpha_2 & \alpha_1 & 0 \end{pmatrix}$. In this representation differential system (8) becomes

$$\text{(9)} \qquad \frac{dK}{dt} \; = \; [\Omega, K] + [A, P], \qquad \frac{dP}{dt} \; = \; [\Omega, P] + k[A, K]$$

in which

$$\widehat{K} \; = \; \begin{pmatrix} H_1 \\ H_2 \\ H_3 \end{pmatrix}, \quad \widehat{\Omega} \; = \; \begin{pmatrix} \frac{1}{c_1} H_1 \\ \frac{1}{c_2} H_2 \\ \frac{1}{c_3} H_3 \end{pmatrix}, \quad \widehat{P} \; = \; \begin{pmatrix} h_1 \\ h_2 \\ h_3 \end{pmatrix}, \quad \text{and} \quad \widehat{A} \; = \; \begin{pmatrix} a_1 \\ a_2 \\ a_3 \end{pmatrix}.$$

The characteristic polynomial of the matrix $U$ in (7) is given by

$$\lambda^4 + \lambda^2(\|\widehat{H}\|^2 + k\|\hat{h}\|^2) + (\widehat{H} \cdot \hat{h})^2 \; ;$$

hence,

$$\text{(10)} \qquad K_2 \; = \; \|\hat{h}\|^2 + k\|\widehat{H}\|^2 \quad \text{and} \quad K_3 \; = \; \hat{h} \cdot \widehat{H}$$

are the Casimir functions on $\mathfrak{g}$. Being constant on each co-adjoint orbit of $G$, they Poisson commute with any function on $\mathfrak{g}^*$, and in particular they Poisson commute with each other. Since $\mathfrak{g}$ is of rank 2, it follows that in addition to $K_2$ and $K_3$ there are two extra integrals of motion for $H$ by our preceding observations about right-invariant vector fields. Together with $H$ these functions constitute five independent integrals of motion, all Poisson commuting with each other. So $H$ will be completely integrable just in case when there is one more independent integral that Poisson commutes with $H$.

*The Euclidean case.* The group of motions $E_3$ is not semisimple, hence the Hamiltonian equations cannot be written in the Lax-pair form as in (6). The following bilinear (but not invariant) form reveals the connections with the equations for the heavy top:

$$\langle A, B \rangle \ = \ \sum_{i=1}^{3} a_i \bar{a}_i + b_i \bar{b}_i,$$

with

$$A \ = \ \sum_{i=1}^{3} a_i A_i + b_i B_i \quad \text{and} \quad B = \sum_{i=1}^{3} \bar{a}_i A_i + \bar{b}_i B_i.$$

Relative to this form every $p = \sum_{i=1}^{n} h_i B_i^* + H_i A_i^*$ in $\mathfrak{g}^*$ is identified with

$$U \ = \ \begin{pmatrix} 0 & 0 & 0 & 0 \\ h_1 & 0 & -H_3 & H_2 \\ h_2 & H_3 & 0 & -H_1 \\ h_3 & -H_2 & H_1 & 0 \end{pmatrix} \quad \text{in} \quad \mathfrak{g}.$$

Then along an extremal curve $(g(t), U(t))$ of $H$, functions $F_L(g, U) = \langle U, g^{-1} L g \rangle$ are constant for each $L$ in $\mathfrak{g}$. Upon expressing $g = \begin{pmatrix} 1 & 0 \\ x & R \end{pmatrix}$ in terms of the translation $x$ and the rotation $R$, we see that $F_L$ becomes a function of the variables $x, R, \hat{h}, \widehat{H}$ and is given by

$$F_L(x, R, \hat{h}, \widehat{H}) \ = \ \hat{h} \cdot (R^{-1}(v + Vx) + \widehat{H} \cdot R^{-1} \widehat{V}$$

with

$$L \ = \ \begin{pmatrix} 0 & 0 & 0 & 0 \\ v_1 & 0 & -V_3 & V_2 \\ v_2 & V_3 & 0 & -V_1 \\ v_3 & -V_2 & V_1 & 0 \end{pmatrix}, \quad \text{and} \quad \widehat{V} \ = \ \begin{pmatrix} \widehat{V_1} \\ \widehat{V_2} \\ \widehat{V_3} \end{pmatrix}.$$

The Lie algebra of $E_3$ is of rank 3, because all translations commute. Taking $L$ in the space of translations amounts to taking $\widehat{V} = 0$, and so functions

$$F_v(x, R, \hat{h}, \widehat{H}) \ = \ \hat{h} \cdot R^{-1} v$$

Poisson commute with each other, and are also integrals of motion for any left-invariant Hamiltonian. The functions $F_v$ form a three dimensional space with $F_1 = R\hat{h} \cdot e_1$, $F_2 = R\hat{h} \cdot e_2$, $F_3 = R\hat{h} \cdot e_3$ a basis for such a space. The elements of this basis are not functionally independent because of the following relation:

$$F_1^2 + F_2^2 + F_3^2 \;=\; h_1^2 + h_2^2 + h_3^2.$$

Consequently, the functions $F_v$ give at most two independent integrals of motion.

By using the Poisson bracket, Table 1 one shows that the differential equations for $U(t)$ can be written as

$$\frac{d\hat{h}(t)}{dt} \;=\; \hat{h}(t) \times \widehat{\Omega}(t), \quad \frac{d\widehat{H}}{dt} \;=\; \widehat{H}(t) \times \widehat{\Omega}(t) + \hat{h}(t) \times a.$$

Hence, the Hamiltonian equations in this case coincide with equations (8) and (9) for $k = 0$. Consequently, the conservation laws defined by (10) apply to $k = 0$ and we get that

$$K_2 \;=\; \|\hat{h}(t)\|^2 \;=\; \text{constant}, \quad \text{and} \quad K_3 \;=\; \widehat{H}(t) \cdot \hat{h}(t) \;=\; \text{constant}$$

along the integral curves of $\vec{H}$.

Together with $H, K_2, K_3$ and any two functions among $F_1, F_2, F_3$ account for five independent integrals of motion all in involution with each other, and the question of complete integrability in this case also reduces to finding one extra integral of motion.

Having shown that the equations (i) in the introduction coincide with the Hamiltonian equations (9) we now turn to the integrable cases. Since the extra integrals of motion for the three cases mentioned in the introduction are evident, we shall go directly to the case discovered by Kowalewski.

So assume that $c_1 = c_2 = 2c_3$ and that $a_3 = 0$. Normalize the constants so that $c = c_2 = 2$ and $c_3 = 1$. Then, equations (8) become:

$$\frac{dh_1}{dt} = H_3 h_2 - \frac{1}{2} H_2 h_3 - k a_2 H_3, \qquad\qquad \frac{dH_1}{dt} = \frac{1}{2} H_2 H_3 - a_2 h_3 \quad,$$

$$\frac{dh_2}{dt} = \frac{1}{2} H_1 h_3 - H_3 h_1 + k a_1 H_3, \qquad\qquad \frac{dH_2}{dt} = -\frac{1}{2} H_1 H_3 + a_1 h_3,$$

$$\frac{dh_3}{dt} = \frac{1}{2} (H_2 h_1 - H_1 h_2) + k(a_2 H_1 - a_1 H_2), \qquad \frac{dH_3}{dt} = a_2 h_1 - a_1 h_2.$$

Set $z \;=\; \frac{1}{2}(H_1 + iH_2), \quad w = h_1 + ih_2$, and $a = a_1 + ia_2$. Then

$$\frac{dz}{dt} \;=\; \frac{-i}{2}(H_3 z - a h_3) \quad \text{and} \quad \frac{dw}{dt} = i(h_3 z - H_3 w + k H_3 a).$$

Let $q = z^2 - a(w - ka)$. Then,

$$\frac{dq}{dt} = 2z\frac{dz}{dt} - a\frac{dw}{dt} = -i(H_3 z^2 - ah_3 z) - ia(h_3 z - H_3 w + kH_3 a)$$
$$= -iH_3(z^2 - aw + ka^2) = -iH_3(t)q(t).$$

Denoting by $\bar{q}$ the complex conjugate of $q$ we get that $\frac{d\bar{q}}{dt} = iH_3(t)\bar{q}(t)$, and hence

$$\frac{d}{dt}q(t)\bar{q}(t) = -iH_3 q\bar{q} + iH_3 q\bar{q} = 0.$$

Thus,

$$q(t)\bar{q}(t) = |q(t)|^2 = \text{constant}.$$

Hence $|z^2 - a(w - ak)|^2$ is the required integral of motion.

We shall refer to this case as the Kowalewski case.

## 2. The Kowalewski case: Reductions and eliminations

It will be convenient to rescale the coordinates so that the constant $a$ is reduced to 1. Let

$$x = \frac{\bar{a}}{|a|^2}z\left(\frac{t}{|a|}\right), \ x_3 = \frac{1}{|a|}H_3\left(\frac{t}{|a|}\right), \ y = \frac{\bar{a}}{|a|^2}w\left(\frac{t}{|a|}\right), \ y_3 = \frac{1}{|a|}h_3\left(\frac{t}{|a|}\right).$$

It follows from the previous page that

$$(11) \qquad \frac{dx}{dt} = -\frac{i}{2}(x_3 x - y_3), \quad \frac{dy}{dt} = i(y_3 x - x_3 y + kx_3) ,$$
$$\frac{dx_3}{dt} = \text{Im}\, y, \quad \text{and} \quad \frac{dy_3}{dt} = (\text{Im}\, x\bar{y} + 2k\text{Im}\, \bar{x})$$

are the extremal equations in our new coordinates.

The integrals of motion in these coordinates become:

$$(12) \qquad H = \frac{1}{4}(H_1^2 + H_2^2) + \frac{1}{2}H_3^2 + a_1 h_1 + a_2 h_2$$
$$= z\bar{z} + \frac{1}{2}H_3^2 + \text{Re}\, aw = |a|^2(x\bar{x} + \frac{1}{2}x_3^2 + \text{Re}\, y) ,$$
$$K_2 = \|\hat{h}\|^2 + k\|\widehat{H}\|^2 = |a|^2(|y|^2 + y_3^2 + k(4|x|^2 + x_3^2)) ,$$
$$K_3 = \hat{h}\cdot\widehat{H} = |a|^2(2\text{Re}\, x\bar{y} + x_3 y_3) ,$$
$$K_4^2 = q\bar{q} = |a|^2|(x^2 - (y - k)|^2$$

and therefore we may assume that $a = 1$. This rescaling reveals that system (11) is invariant under the involution

$$\sigma(x, y, x_3, y_3) = (\bar{x}, \bar{y}, -x_3, -y_3).$$

We shall now assume that the constants $H, K_2, K_3, K_4$ are fixed, and use $V$ to denote the manifold defined by equations (12). Now $V$ is a two dimensional real variety, contained in $\mathbb{R}^6$, that can be conveniently parametrized by one complex variable according to the following theorem.

THEOREM 1. $V$ *is contained in the set of all complex numbers* $x$ *and* $q$ *that satisfy*

$$(13) \qquad P(x)\bar{q} + P(\bar{x})q + R_1(x, \bar{x}) + K_4^2(x - \bar{x})^2 = 0$$

*with*

$$P(x) = \widetilde{K}_2 - 2K_3 x + 2Hx^2 - x^4, \qquad and$$

$$(14)$$
$$R_1(x, \bar{x}) = (\widetilde{H}\widetilde{K}_2 - K_3^2) + 2K_3 k(x + \bar{x}) + (2\widetilde{H}k - 3K_2)(x^2 + \bar{x}^2)$$
$$+ 2K_3 x\bar{x}(x + \bar{x}) - \widetilde{H}x^2\bar{x}^2 + (\widetilde{H}k - 2K_2)(x - \bar{x})^2 ,$$

*where* $\widetilde{H} = 2H - 2k$, $\widetilde{K}_2 = K_2 - k\widetilde{H} - K_4^2$.

*Proof.* Equation (13) is a consequence of eliminating $x_3$ and $y_3$ from the constraints (12) as follows:

Begin by expressing the integrals of motion in terms of $x, \bar{x}, q, \bar{q}, x_3$ and $y_3$. Putting $y = x^2 - q + k$, in the expression for $H$ leads to

$$H = x\bar{x} + \frac{1}{2}x_3^2 + \text{Re}(x^2 - q + k) = x\bar{x} + \frac{1}{2}x_3^2 + \frac{1}{2}(x^2 + \bar{x}^2 - (q + \bar{q}) + 2k).$$

This relation simplifies to $2H - 2k = \widetilde{H} = (x + \bar{x})^2 - (q + \bar{q}) + x_3^2$.

Then, $K_2 = x^2\bar{x}^2 + k^2 + y_3^2 + k\widetilde{H} + 2kx\bar{x} - (x^2\bar{q} + \bar{x}^2 q)$ and hence,

$$K_2 - k\widetilde{H} - k^2 = \widetilde{K}_2 = x^2\bar{x}^2 + y_3^2 + 2kx\bar{x} - (x^2\bar{q} + \bar{x}^2 q).$$

Finally, $K_3 = (x\bar{x} + k)(x + \bar{x}) - (x\bar{q} + \bar{x}q) + x_3 y_3$.

Eliminating $y_3$ and $x_3$ from the preceding relations leads to:

$$(14a)$$
$$(K_3 - (x\bar{x} + k)(x + \bar{x}) + (x\bar{q} + \bar{x}q))^2 = x_3^2 y_3^2$$
$$= (\widetilde{K}_2 + \bar{x}^2 q + x^2\bar{q} - x^2\bar{x}^2 - 2kx\bar{x})(\widetilde{H} - (q + \bar{q}) - (x + \bar{x})^2).$$

The homogeneous terms of degree two in $q$ and $\bar{q}$ in the preceding expression reduce to

$$(x\bar{q} + \bar{x}q)^2 - (\bar{x}^2 q + x^2\bar{q})(q + \bar{q})$$
$$= x^2\bar{q}^2 + \bar{x}^2 q^2 + 2x^2\bar{x}^2 q\bar{q} - (\bar{x}^2 q^2 + x^2 q\bar{q} + \bar{x}^2 q\bar{q} + x^2\bar{q}^2)$$
$$= -K_4^2(x - \bar{x})^2.$$

Therefore, relation (14a) can be reduced to

$$P\bar{q} + \overline{P}q + \widehat{R} = 0$$

for suitable polynomials $\widehat{R}$ and $P$ in the variables $x$ and $\bar{x}$. It follows that

$$
\begin{aligned}
P &= (\widetilde{K}_2 - x^2\bar{x}^2 - 2kx\bar{x}) + x^2(\widetilde{H} - (x+\bar{x})^2) - 2(K_3 - (x+\bar{x})(x\bar{x}+k))x \\
&= \widetilde{K}_2 - 2K_3x + x^2(\widetilde{H} - 2k) - x^4 - x^2\bar{x}^2 - 2kx\bar{x} - x^2\bar{x}^2 - 2x^3\bar{x} \\
&\quad + 2kx\bar{x} + 2x^2\bar{x}^2 + 2x^3\bar{x} \\
&= \widetilde{K}_2 - 2K_3x + 2Hx^2 - x^4.
\end{aligned}
$$

Thus $P$ is a polynomial of degree 4 in the variable $x$ only.

Then, $\widehat{R}(x,\bar{x}) = R_1(x,\bar{x}) + K_4^2(x-\bar{x})^2$ with

$$
R_1 = (\widetilde{H} - (x+\bar{x})^2)(\widetilde{K}_2 - x^2\bar{x}^2 - 2kx\bar{x}) - (K_3 - (x+\bar{x})(x\bar{x}+k))^2.
$$

The expression for $R_1$ further simplifies to

$$
\begin{aligned}
R_1 &= (\widetilde{H}\widetilde{K}_2 - K_3^2) + 2kK_3(x+\bar{x}) + (2\widetilde{H}k - 3K_2)(x^2+\bar{x}^2) + 2K_3x\bar{x}(x+\bar{x}) \\
&\quad - \widetilde{H}^2x^2\bar{x}^2 + (\widetilde{H}k - 2K_2)(x-\bar{x})^2
\end{aligned}
$$

by a straightforward calculation.

Equation (13) identifies $x$ as the pivotal variable, in terms of which the extremal equations can be integrated by quadrature. For then $q$ is the solution of (13), and the remaining variables are given by

$$
x_3^2 = \widetilde{H} - (x+\bar{x})^2 + (q+\bar{q}), \quad \text{and} \quad y_3^2 = \widetilde{K}_2 - x^2\bar{x}^2 + x^2q + \bar{x}^2q - 2kx\bar{x}.
$$

THEOREM 2. *Each extremal curve $x(t)$ satisfies the following differential equation*:

$$
\tag{15} -4\left(\frac{dx}{dt}\right)^2 = P(x) + q(t)(x-\bar{x})^2,
$$

*with $P(x)$ as in the previous theorem.*

*Proof.* It follows from equation (11) that $-4\left(\frac{dx}{dt}\right)^2 = (x_3x - y_3)^2$.

$$
\begin{aligned}
(x_3x - y_3)^2 &= x_3^2x^2 - 2x_3y_3x + y_3^2 = (\widetilde{H} + (q+\bar{q}) - (x+\bar{x})^2)x^2 \\
&\quad - 2x(K_3 - (x+\bar{x})(x\bar{x}+k) + (x\bar{q} + \bar{x}q)) \\
&\quad + (\widetilde{K}_2 - 2kx\bar{x} - x^2\bar{x}^2 + x^2\bar{q} + \bar{x}^2q) \\
&= \widetilde{K}_2 - 2K_3x + \widetilde{H}x^2 + (x^2(q+\bar{q}) - (x+\bar{x})^2x^2 \\
&\quad + 2x(x+\bar{x})(x\bar{x}+k) \\
&\quad - 2x(x\bar{q} + \bar{x}q) - 2kx\bar{x} - x^2\bar{x}^2 + x^2\bar{q} + \bar{x}^2q) \ .
\end{aligned}
$$

But then,

$$
\begin{aligned}
&x^2(q+\bar{q}) - (x+\bar{x})^2x^2 + 2x(x+\bar{x})(x\bar{x}+k) - 2x(x\bar{q} + \bar{x}q) \\
&\quad - 2kx\bar{x} - x^2\bar{x}^2 + x^2\bar{q} + \bar{x}^2q = q(x-\bar{x})^2 - x^4 + 2kx^2,
\end{aligned}
$$

and therefore

$$(x_3 x - y_3)^2 = \widetilde{K}_2 - 2K_3 x + \widetilde{H} x^2 + 2kx^2 - x^4 + q(x - \bar{x})^2$$
$$= P(x) + q(x - \bar{x})^2.$$

THEOREM 3. *Let* $R_0(x, \bar{x}) = \widetilde{K}_2 - K_3(x + \bar{x}) + H(x^2 + \bar{x}^2) - x^2\bar{x}^2 - (H - k)(x - \bar{x})^2$. *Then*

(16)    $R_0(x, x) = P(x), \quad and \quad R_0^2(x, \bar{x}) + (x - \bar{x})^2 R_1(x, \bar{x}) = P(x)P(\bar{x})$

*where* $R_1$ *has the same meaning as in Theorem 1.*

*Proof.* Let $\zeta = x_3 x - y_3$. We shall first show that $R_0 = \zeta\bar{\zeta}$.

$$\zeta\bar{\zeta} = (x_3 x - y_3)(x_3\bar{x} - y_3) = x_3^2 x\bar{x} - x_3 y_3(\bar{x} + x) + y_3^2$$
$$= (\widetilde{H} + (q + \bar{q}) - (x + \bar{x})^2)x\bar{x} + (\widetilde{K}_2 - x^2\bar{x}^2 - 2kx\bar{x} + \bar{x}^2 q + x^2\bar{q})$$
$$- (K_3 + x\bar{q} + \bar{x}q - (x + \bar{x})(x\bar{x} + k))(x + \bar{x})$$
$$= \widetilde{K}_2 - K_3(x + \bar{x}) + k(x + \bar{x})^2 - 2kx\bar{x} + \widetilde{H}x\bar{x} - x^2\bar{x}^2$$
$$+ (x\bar{x}(q + \bar{q}) - x\bar{x}(x + \bar{x})^2 + \bar{x}^2 q + x^2\bar{q}$$
$$- (x\bar{q} + \bar{x}q)(x + \bar{x}) + x\bar{x}(x + \bar{x})^2).$$

The above expression reduces to

$$\widetilde{K}_2 - K_3(x + \bar{x}) + k(x^2 + \bar{x}^2) + \widetilde{H}x\bar{x} - x^2\bar{x}^2$$

because

$$x\bar{x}(q + \bar{q}) + \bar{x}^2 q + x^2\bar{q} - (x\bar{q} + \bar{x}q)(x + x) = x\bar{x}(q + \bar{q}) - x\bar{x}(q + \bar{q}) = 0.$$

Since $2k + \widetilde{H} = 2H$, the preceding expression can also be written as

$$\widetilde{K}_2 - K_3(x + \bar{x}) + H(x^2 + y^2) - x^2\bar{x}^2 - (H - k)(x - \bar{x})^2$$

showing that $R_0 = \zeta\bar{\zeta}$.

It follows from the proof of Theorem 1 that $\zeta^2 = P(x) + q(x - \bar{x})^2$. Therefore, $R_0(x, x) = P(x)$ and,

$$R_0^2(x, \bar{x}) = (\zeta\bar{\zeta})^2 = \zeta^2\bar{\zeta}^2 = (P(x) + q(x - x)^2)(P(\bar{x}) + \bar{q}(x - \bar{x})^2)$$
$$= P(x)P(\bar{x}) + P(x)\bar{q}(x - \bar{x})^2 + P(\bar{x})q(x - \bar{x})^2 + q\bar{q}(x - \bar{x})^4.$$

Thus

$$R_0^2(x, \bar{x}) = P(x)P(\bar{x}) - R_1(x, \bar{x})(x - \bar{x})^2,$$

because $P(x)\bar{q} + P(\bar{x})q = -R_1(x, \bar{x}) - K_4^2(x - \bar{x})^2$ as can be seen from relations (13) in Theorem 1. Our theorem is proved.

The essential relations obtained by the preceding calculations, can now be summarized for further reference as follows.

For each choice of constants of motion, the algebraic variety defined by equations (12) is parametrized by a single complex variable $x$ through the relations

$$P(\bar{x})q^2 + (R_1(x,\bar{x}) + K_4^2(x-\bar{x})^2)q + K_4^2 P(x) \ = \ 0 \quad \text{and} \quad q\bar{q} = K_4^2,$$

with $P(x)$ a polynomial of degree four, and $R_1(x,\bar{x})$ the form defined by (14). The form $R_1$ has a companion form $R_0(x,y)$, also of degree four, that satisfies

$$R_0(x,x) \ = \ P(x) \quad \text{and} \quad R_0^2(x,\bar{x}) + (x-\bar{x})^2 R_1(x,\bar{x}) \ = \ P(x)P(\bar{x}).$$

Finally, $-4\left(\frac{dx}{dt}\right)^2 \ = \ P(x) = q(x-\bar{x})^2$ is the extremal differential equation that needs to be solved.

Apart from the more general nature of the constants that occur in the foregoing expressions, all of the above relations are the same as in the original paper of Kowalewski from 1889. It will be convenient to refer to the above relations as Kowalewski's relations.

## 3. Lemniscatic integrals and addition formulas of L. Euler and A. Weil

We shall now show that Kowalewski's relations provide remarkable insights into the theory of elliptic curves and elliptic integrals starting with the very beginning of the subject with the work of G.S. Fagnano in 1718 concerning the arc of a lemniscate, and the subsequent extensions of Euler concerning the solutions of a more general differential equation

$$\frac{dx}{dt} \ = \ \sqrt{P(x)}$$

where $P(x)$ is an arbitrary fourth degree polynomial.([12])

To begin with, note that when the constant of motion $K_4$ is equal to zero, then

$$-4\left(\frac{dx}{dt}\right)^2 \ = \ P(x).$$

Thus our extremal equations already contain the cases studied by Euler in connection with the length of an arc of a lemniscate. It is therefore not surprising that our study must be closely related to the work of Euler. However, it is remarkable that the fundamental relation

$$R_0^2(x,\bar{x}) + (x-\bar{x})^2 R_1(x,\bar{x}) \ = \ P(x)P(\bar{x})$$

of Kowalewski provides an easy access to Euler's famous results concerning the solutions of

$$\frac{dx}{\sqrt{P(x)}} \pm \frac{dy}{\sqrt{P(y)}} = 0.$$

To explain Euler's results and its relevance to this paper we shall assume that

$$P(x) = A + 4Bx + 6Cx^2 + 4Dx^3 + Ex^4$$

is a general polynomial of degree four. We shall also denote by $R(x, y) = A + 2B(x + y) + 3C(x^2 + y^2) + 2Dxy(x + y) + Ex^2y^2$ a particular form that satisfies $R(x, x) = P(x)$.

THEOREM 4.

$$\begin{aligned}
\widehat{R}(x, y) =\ & -4B^2 + 4(AD - 3BC)(x + y) \\
& + 2(AE + 2BD - 9C^2)(x^2 + y^2) + 4(BE - 3CD)xy(x + y) \\
& - 4D^2x^2y^2 - (AE + 4BD - 9C^2)(x - y)^2
\end{aligned}$$

is the unique form that satisfies

$$R^2(x, y) + (x - y)^2 \widehat{R}(x, y) = P(x)P(y) .$$

In particular, $\widehat{R}(x, x)$ defines a polynomial $Q$ given by

$$\begin{aligned}
Q(x) = \ & 4(-B^2 + 2(AD - 3BC)x^2 + (AE + 2BD - 9C^2)x^2 \\
& + 2(BE - 3CD)x^3 - D^2x^4).
\end{aligned}$$

*Proof.* Consider $F(x, y) = P(x)P(y) - R^2(x, y)$. Because $P(x) = R(x, x)$, $P'(x) = 2\frac{\partial R}{\partial x}(x, x)$. Therefore, $\frac{\partial F}{\partial x}(x, y) = P'(x)P(y) - 2R(x, y)\frac{\partial R}{\partial x}(x, y)$, and for $x = y$

$$\frac{\partial F}{\partial x}(x, x) = P'(x)P(x) - 2P(x)\frac{1}{2}P'(x) = 0.$$

It follows by an analogous argument that $\frac{\partial F}{\partial y}(x, y) = 0$ for $y = x$. Thus both $F$ and its partial derivations $\frac{\partial F}{\partial x}$ and $\frac{\partial F}{\partial y}$ vanish at $x = y$. Consequently, $F(x, y) = (x - y)^2 \widehat{R}(x, y)$ for some binary form $\widehat{R}(x, y)$.

Now, by differentiation of $R^2(x, y) + (x - y)^2 \widehat{R}(x, y) = P(x)P(y)$,

$$2R(x, y)\frac{\partial R}{\partial x}(x, y) + 2(x - y)\widehat{R}(x, y) + (x - y)^2\frac{\partial \widehat{R}}{\partial x}(x, y) = P'(x)P(y),$$

and

$$\begin{aligned}
2\left(\frac{\partial R}{\partial y}\frac{\partial R}{\partial x} + R\frac{\partial^2 R}{\partial y \partial x}\right) - 2\widehat{R} + 2(x - y)\frac{\partial \widehat{R}}{\partial y} - 2(x - y)\frac{\partial \widehat{R}}{\partial x} \\
+ (x - y)^2\frac{\partial^2 \widehat{R}}{\partial y \partial x} = P'(x)P'(y).
\end{aligned}$$

Let $Q(x)$ denote $\widehat{R}(x,x)$, and set $x = y$ in the above expression. It follows that

$$2\left(\left(\frac{\partial R}{\partial x}(x,x)\right)^2 + P(x)\frac{\partial^2 R}{\partial y \partial x}(x,x)\right) - 2Q(x) = P'(x)^2,$$

which upon substituting $2\frac{\partial R}{\partial x}(x,x) = P'(x)$ reduces to

$$\frac{1}{2}P'(x)^2 + 2P(x)\frac{\partial^2 R}{\partial y \partial x}(x,x) - 2Q(x) = P'(x)^2.$$

Solving for $Q$ gives

$$\begin{aligned}
Q(x) &= P(x)\frac{\partial^2 R}{\partial y \partial x}(x,x) - \frac{1}{4}P'(x)^2 \\
&= (A + 4Bx + 6Cx^2 + 4Dx^2 + Ex^4)(8Dx + 4Ex^2) \\
&\quad - \frac{1}{4}(4B + 12Cx + 12Dx^2 + 4Ex^3)^2 \\
&= -4B^2 + (8AD - 24BC)x + (32BD + 4AE - 24BD - 36C^2)x^2 \\
&\quad + (48CD + 12BE - 72CD - 8BE)x^3 \\
&\quad + (32D^2 + 24CE - 36D^2 - 24EC)x^4 \\
&= -4B^2 + 2(4AD - 12BC)x + (8BD + 4AE - 36C^2)x^2 \\
&\quad + 2(4BE - 12CD)x^3 - 4D^2x^4 \\
&= A' + 4B'x + 6C'x^2 + 4D'x^3 + E'x^4.
\end{aligned}$$

Therefore,

$$\widehat{R}(x,y) = A' + 2B'(x+y) + \gamma(x^2+y^2) + 2\delta xy + 2D'xy(x+y) + E'x^2y^2$$

with $2\gamma + 2\delta = 6C'$, where $\gamma$ and $\delta$ are to be determined from the fundamental relation $R^2 + (x-y)^2\widehat{R} = P(x)P(y)$. Upon equating the homogeneous terms of degree 4 in the fundamental relation we get:

$$9C^2(x^2+y^2)^2 + 8BDxy(x+y)^2 + 2AEx^2y^2 - \gamma(x-y)^2(x^2+y^2) + 2\delta xy(x-y)^2$$
$$= AEy^4 + 16BDxy^3 + 36C^2x^2y^2 + 16BDx^3y + AEx^4.$$

It follows that

$$AE(x+y)^2 - 9C^2(x+y)^2 + 8BDxy = \gamma(x^2+y^2) + 2\delta xy$$

and therefore

$$\gamma = AE - 9C^2, \quad \text{and} \quad \delta = AE + 4BD - 9C^2.$$

The substitution of these values in the expression for $\widehat{R}$ leads to

$$\begin{aligned}
\widehat{R}(x,y) \;=\; & -4B^2 + (4AD - 12BC)(x+y) + (4BD + 2AE - 18C^2)(x^2 + y^2) \\
& + (4BE - 12CD)xy(x+y) - 4D^2x^2y^2 \\
& - (AE + 4BD - 9C^2)(x-y)^2
\end{aligned}$$

and the proof of the theorem is finished.

Let $R_\theta(x,y) = R(x,y) - \theta(x-y)^2$ denote the most general biquadratic form that satisfies $R_\theta(x,x) = P(x)$ with $\theta$ an arbitrary parameter. Then $\Phi_\theta(x,y) = -(x-y)^2\theta^2 + 2R(x,y)\theta + \widehat{R}(x,y)$ is the unique form that satisfies

$$(17) \qquad R_\theta^2(x,y) + (x-y)^2\Phi_\theta(x,y) \;=\; P(x)P(y).$$

*Remark.* For the relations of Kowalewski obtained earlier in the paper, $P = \widetilde{K}_2 - 2K_3x + 2Hx^2 - x^4$, and the forms $R_0$ and $R_1$ are given by $R_0 = R_\theta$, with $\theta = H - k$, and $R_1 = -(x-y)^2(H-k)^2 + 2R(x,y)(H-k) + \widehat{R}$.

Now, to return to the general case, $\Phi_\theta$ is symmetric with respect to $x$ and $y$ and can be written as

$$\Phi_\theta(x,y) = a_\theta(x)y^2 + 2b_\theta(x)y + c_\theta(x) \;=\; a_\theta(y)x^2 + 2b_\theta(y)x + c_\theta(y)$$

for some quadratic expressions $a_\theta, b_\theta$ and $c_\theta$.

Writing $\widehat{R} = \alpha + 2\beta(x+y) + 3\gamma(x^2+y^2) + 2\delta xy(x+y) + \epsilon x^2 y^2 - \zeta(x-y)^2$ we have

$$\begin{aligned}
\Phi_\theta \;=\; & -(x-y)^2(\theta^2 + \zeta) + 2(A + 2B(x+y) + 3C(x^2+y^2) \\
& + 2Dxy(x+y) + Ex^2y^2)\theta + (\widehat{R} + \zeta(x-y)^2),
\end{aligned}$$

and therefore,

$$\begin{aligned}
a_\theta(x) \;&=\; (2E\theta + \epsilon)x^2 + (4D\theta + 2\delta)x + (6C\theta + 3\gamma - (\theta^2 + \zeta)) , \\
c_\theta(x) \;&=\; (6C\theta + 3\gamma - (\theta^2 + \zeta))x^2 + (4B\theta + 2\beta)x + (2A\theta + \alpha) , \\
b_\theta(x) \;&=\; (2D\theta + \delta)x^2 + (\theta^2 + \zeta)x + (2B\theta + \beta) .
\end{aligned}$$

After substitution of the values for $\alpha, \beta, \gamma, \delta, \epsilon, \zeta$ given by Theorem 4, the preceding expressions become

$$\begin{aligned}
a_\theta(x) \;&=\; (2E\theta - 4D^2)x^2 + (4D\theta + 4(BE - 3CD))x + AE - (\theta - 3C)^2 , \\
c_\theta(x) \;&=\; (AE - (\theta - 3C)^2)x^2 + (4B\theta + 4(AD - 3BC))x + (2A\theta - 4B^2) , \\
b_\theta(x) \;&=\; (2D\theta + 2(BE - 3CD))x^2 + (\theta^2 - 9C^2 + AE + 4BD)x \\
& \quad + (2B\theta + 2AD - 3BC).
\end{aligned}$$

THEOREM 5. *Let $G_\theta$ denote the discriminant $b_\theta^2 - a_\theta c_\theta$. Then $G_\theta(x) = p(\theta)P(x)$ with*

$$(18) \qquad p(\theta) \;=\; 2\theta(\theta - 3C)^2 + 2\theta(4BD - AE) + 4B^2E + 4AD^2 - 24BCD.$$

*Proof.* Let $x$ be any root of $P(x)$. Then,

$$(R(x,y) - \theta(x-y)^2)^2 + (x-y)^2(-(x-y)^2\theta^2 + 2R(x,y)\theta + \widehat{R}(x,y)) = 0.$$

When

$$\Phi_\theta(x,y) = 0, \qquad R(x,y) - \theta(x-y)^2 = 0$$

and therefore $\Phi_\theta(x,y) = 0$ has a double root. Therefore $G_\theta = 0$. This argument shows that $P(x)$ is a factor of $G_\theta$. Since $G_\theta$ is a polynomial of degree four in $x$, it follows that $G_\theta(x) = p(\theta)P(x)$ for some polynomial $p(\theta)$. It now follows by an easy calculation, using the explicit expressions of $a_\theta, b_\theta$ and $c_\theta$, that $p(\theta)$ is given by expression (18). This ends the proof.

Expression (18) is most naturally linked with the cubic elliptic curve $\Gamma = \{(\xi, \eta) : \eta^2 = 4\xi^3 - g_2\xi - g_3\}$, where $g_2$ and $g_2$ are the invariants of $\mathcal{C} = \{(x, u) : u^2 = P(x)\}$ explicitly given as follows:

$$(19) \quad g_2 = AE - 4BD + 3C^2, \quad g_3 = ACE + 2BCD - AD^2 - B^2E - C^3.$$

The identification is obtained as follows: Let $\theta = 2(\xi + C)$. Then, $p(\xi) = 4(4\xi^3 - (AE - 4BD + 3C^2)\xi - (ACE + 2BCD - AD^2 - B^2E - C^3))$.

Letting $\eta^2 = \frac{p}{4}$ we obtain $\eta^2 = 4\xi^3 - g_2\xi - g_3$.

The next theorem is a paraphrase of the classical results of Euler ([3]).

THEOREM 6. $\Phi_\theta(x,y) = 0$ *is a solution for* $\frac{dx}{\sqrt{P(x)}} \pm \frac{dy}{\sqrt{P(y)}} = 0$ *for each number* $\theta$. *Conversely, for every solution* $y(x)$ *of either differential equation* $\frac{dx}{\sqrt{P(x)}} \pm \frac{dy}{\sqrt{P(y)}} = 0$ *there exists a number* $\theta$ *such that* $\Phi_\theta(x,y) = 0$ *is equal to* $y(x)$.

*Proof.* Consider $\Phi_\theta(x,y) = 0$ for an arbitrary number $\theta$. Then

$$(20) \qquad x = \frac{-b_\theta(y) + \rho\sqrt{G_\theta(y)}}{a_\theta(y)} \quad \text{and} \quad y = \frac{-b_\theta(x) + \sigma\sqrt{G_\theta(x)}}{a_\theta(x)}$$

with $G_\theta = b_\theta^2 - a_\theta c_\theta$, and $\rho = \pm 1$, $\sigma = \pm 1$. Hence

$$(21) \qquad \frac{1}{2}\frac{\partial\Phi_\theta}{\partial x} = xa_\theta(y) + b_\theta(y) = \rho\sqrt{G_\theta(y)} \quad \text{and}$$

$$\frac{1}{2}\frac{\partial\Phi_\theta}{\partial y} = ya_\theta(x) + b_\theta(x) = \sigma\sqrt{G_\theta(x)} \ .$$

Using the results of Theorem 5 we get that the curve $y(x)$ of (20) satisfies:

$$\frac{dy}{dx} = -\frac{\partial\Phi_\theta}{\partial x}\bigg/\frac{\partial\Phi_\theta}{\partial y} = \frac{-\rho\sqrt{G_\theta(y)}}{\sigma\sqrt{G_\theta(x)}} = \frac{-\rho\sqrt{p(\theta)P(y)}}{\sigma\sqrt{p(\theta P(x)}},$$

and therefore

$$\frac{dy}{\sqrt{P(y)}} + \frac{\rho}{\sigma}\frac{dx}{\sqrt{P(x)}} = 0.$$

Thus, $y(x)$ is a solution of $\frac{dx}{\sqrt{P(x)}} + \frac{dy}{\sqrt{P(y)}} = 0$ when $\frac{\rho}{\sigma} = 1$. Otherwise, it is a solution of $\frac{dx}{\sqrt{P(x)}} - \frac{dy}{\sqrt{P(y)}} = 0$ .

Suppose now that $y(x)$ is a particular solution of either $\frac{dx}{\sqrt{P(x)}} + \frac{dy}{\sqrt{P(y)}} = 0$ or $\frac{dx}{\sqrt{P(x)}} - \frac{dy}{\sqrt{P(y)}} = 0$. Let $x = a$, $y = b$ be any point such that $b = y(a)$. We need to show that there exists $\theta$ such that $\Phi_\theta(a, b) = 0$, or that

$$\theta = \frac{R(a,b) \pm \sqrt{R^2(a,b) + (a-b)^2 \widehat{R}(a,b)}}{(a-b)^2}.$$

Since

$$R^2(a,b) + (a-b)^2 \widehat{R}(a,b) = P(a)P(b)$$

we get that

$$(22) \qquad \theta = \frac{R(a,b) \pm \sqrt{P(a)P(b)}}{(b-a)^2}.$$

The appropriate sign for $\theta$ depends on whether $y(x)$ is a solution of $\frac{dx}{\sqrt{P(x)}} + \frac{dy}{\sqrt{P(y)}} = 0$, or $\frac{dx}{\sqrt{P(x)}} - \frac{dy}{\sqrt{P(y)}}$. The correct way of choosing the sign for $\theta$ will be made clear through the discussion of the related addition formulas of A. Weil.

*Addition formulas of A. Weil.* A. Weil points out in his paper [13] that the results of Euler (Theorem 6) have algebraic interpretations that may be used to define a group structure on $\Gamma \cup \mathcal{C}$. Recall that

$$\Gamma = \{(\xi, \eta) : \eta^2 = 4\xi^3 - g_2\xi - g_3\} \quad \text{and} \quad \mathcal{C} = \{(x, u) : u^2 = P(x)\}.$$

It turns out, quite remarkably, that these algebraic observations of Weil provide easy explanations for some of the formulas used by Kowalewski, and for that reason it will be necessary to explain Weil's interpretation of Euler's results.

With each solution $\Phi_\theta(x, y) = 0$ Weil associates two transformations from $\mathcal{C}$ into $\mathcal{C}$, depending whether they change $\frac{dx}{u}$ into $\frac{dy}{v}$, or into $-\frac{dy}{v}$, each parametrized by the points of $\Gamma$ rather than $\theta$. More precisely, let $P = (\xi, \eta)$ be any point of $\Gamma$. Let $\Phi_\theta$ denote the form corresponding to $\theta = 2(\xi + C)$. For each point $M = (x, u)$ of $\mathcal{C}$ define two points $N_1 = (y_1, v_1)$ and $N_2 = (y_2, v_2)$ on $\mathcal{C}$ by the following formulas:

$$(23(a)) \qquad y_1 = \frac{-b_\theta(x) + 2\eta u}{a_\theta(x)} , \qquad v_1 = -\frac{1}{2\eta}(xa_\theta(y_1) + b_\theta(y_1)) ,$$

$$(23(b)) \qquad y_2 = \frac{-b_\theta(x) - 2\eta u}{a_\theta(x)} , \qquad v_2 = -\frac{1}{2\eta}(xa_\theta(y_2) + b_\theta(y_2)) .$$

It follows that

$$\frac{dy_1}{dx} = -\frac{\partial \Phi_\theta}{\partial x} \Big/ \frac{\partial \Phi_\theta}{\partial y} = \frac{2\eta v_1}{2\eta u}, \quad \text{and} \quad \frac{dy_2}{dx} = -\frac{\partial \Phi_\theta}{\partial x} \Big/ \frac{\partial \Phi_\theta}{\partial y} = \frac{-2\eta v_2}{2\eta u} \ ;$$

consequently, the mapping $(P, M) \to N_1$ takes $\frac{dx}{u}$ into $\frac{dy}{v}$, while the mapping $(P, M) \longrightarrow N_2$ takes $\frac{dx}{u}$ into $-\frac{dy}{v}$.

Following Weil we shall write $P + M = N_1$ and $P - M = N_2$, with the understanding that $-P = (\xi, -\eta)$ and $-M = (x, -u)$. We assume that on $\Gamma$ the group law coincides with the usual group law, with the point at infinity on $\Gamma$ acting as the group identity.

Formula (22) may be used to show that for each $M$ and $N$ on $\mathcal{C}$ there exist points $P$ and $P'$ on $\Gamma$ such that

$$P + M = N \quad \text{and} \quad P' - M = N$$

by the following simple argument: let $P = (\xi, \eta)$ and $P' = (\xi', \eta')$. Note that $P' = 2M$ when $N = M$, while $P$ is the point at infinity when $N = M$. Thus, $\xi = 2(\theta + C)$ and $\xi' = 2(\theta' + C)$ with $\theta$ and $\theta'$ appropriately chosen among

$$\frac{R(x, y) \pm uv}{(y - x)^2}.$$

For $uv = \sqrt{P(x)P(y)}$, $\frac{R(x,y)-uv}{(x-y)^2}$ gives a finite value when $x = y$ because

$$\begin{aligned}
\frac{R(x, y) - uv}{(x - y)^2} &= \frac{R(x, y) - \sqrt{P(x)P(y)}}{(x - y)^2} \\
&= \frac{(R(x, y) - \sqrt{P(x)P(y)})(R(x, y) + \sqrt{P(x)P(y)})}{(x - y)^2 (R(x, y) + \sqrt{P(x)P(y)})} \\
&= \frac{R^2(x, y) - P(x)P(y)}{(x - y)^2 (R(x, y) + \sqrt{P(x)P(y)})} = \frac{-\widehat{R}(x, y)}{R(x, y) + \sqrt{P(x)P(y)}},
\end{aligned}$$

and the latter expression evidently tends to $\frac{-\widehat{R}(x,x)}{2P(x)}$ when $y$ tends to $x$. Thus,

$$\theta' = \frac{R(x, y) - uv}{(x - y)^2}.$$

Then use formula (23(b)) to define $\eta' = -\frac{1}{2v}(xa_{\theta'}(y) + b_{\theta'}(y))$.

The other choice of sign $uv = -\sqrt{P(x)P(y)}$ leads to $\theta' = \frac{R(x,y)+uv}{(x-y)^2}$ by an analogous argument. In such a case define $\eta' = \frac{1}{2v}(xa_{\theta'}(y) + b_\theta(y))$.

The values of $\eta$ corresponding to $\theta$ are given by

$$\eta = -\frac{1}{2v}(xa_\theta(y) + b_\theta(y)) \quad \text{when} \quad uv = \sqrt{P(x)P(y)}, \quad \text{and}$$

$$\eta = \frac{1}{2v}(xa_\theta(y) + b_\theta(y)) \quad \text{when} \quad uv = -\sqrt{P(x)P(y)}.$$

To these addition formulas of Weil we can add their infinitesimal versions that appeared in the paper of Kowalewski.

THEOREM 7. *Let $O = (\xi, \eta)$ and $O' = (\xi', \eta')$ denote points of $\Gamma$, and let $M = (x, u)$ and $N = (y, v)$ denote points of $\mathcal{C}$ related by the formulas $O = N - M, \quad O' = N + M$. Then,*

$$\frac{d\xi'}{\eta'} = \frac{dx}{u} + \frac{dy}{v}, \quad and \quad \frac{d\xi}{\eta} = -\frac{dx}{u} + \frac{dy}{v}.$$

*Proof.* Since $\Phi_\theta(x, y) = -(x - y)^2 \theta^2 + 2R(x, y)\theta + \widehat{R}(x, y)$,

$$d\Phi_\theta(x, y) = (-2(x - y)^2 \theta + 2R)d\theta + \frac{\partial \Phi_\theta}{\partial x}dx + \frac{\partial \Phi_\theta}{\partial y}dy.$$

For $(x, y, \theta)$ for which $\Phi_\theta(x, y) = 0$, $d\Phi_\theta(x, y) = 0$, and $\theta = \frac{R(x,y) \pm uv}{(x-y)^2}$ depending whether $uv = \sqrt{P(x)P(y)}$, or $uv = -\sqrt{P(x)P(y)}$.

Assuming that $uv = \sqrt{P(x)P(y)}$, we have $\theta = \frac{R(x,y) + uv}{(x-y)^2}$ and $\frac{1}{2}\frac{\partial \Phi_\theta}{\partial x} = xa_\theta(y) + b_\theta(y) = -2\eta v$ and $\frac{1}{2}\frac{\partial \Phi_\theta}{\partial y} = ya_\theta(x) + b_\theta(x) = 2\eta u$. Thus,

$$-2uvd\theta - 4\eta v dx + 4\eta u dy = 0,$$

and since $d\theta = 2d\xi$,

$$\frac{d\xi}{\eta} = -\frac{dx}{u} + \frac{dy}{v}.$$

In the remaining case $uv = -\sqrt{P(x)P(y)}$, and $\frac{1}{2}\frac{\partial \phi_\theta}{\partial x} = 2\eta v$, and $\frac{1}{2}\frac{\partial \phi_\theta}{\partial y} = 2\eta u$. Thus, again

$$\frac{d\xi}{\eta} = -\frac{dx}{u} + \frac{dy}{v}.$$

The corresponding differential form for $\xi'$ is obtained by analogous arguments and its derivation will be omitted.

With these formulas behind us we finally come to the integration of the extremal differential equations.

## 4. Kowalewski's integration procedure

In her original paper of 1889 Kowalewski made the following change of variables:

$$s_1 = \frac{R_0(x_1, x_2) - \sqrt{P(x_1)P(x_2)}}{2(x_1 - x_2)^2} + \frac{1}{2}\ell_1$$

and

$$s_2 = \frac{R_0(x_1, x_2) + \sqrt{P(x_1)P(x_2)}}{2(x_1 - x_2)^2} + \frac{1}{2}\ell_1$$

and then claimed that

$$\frac{ds_1}{S(s_1)} = \frac{dx_1}{\sqrt{P(x_1)}} + \frac{dx_2}{\sqrt{P(x_2)}}, \quad \text{and} \quad \frac{ds_2}{S(s_2)} = -\frac{dx_1}{\sqrt{P(x_1)}} + \frac{dx_2}{\sqrt{P(x_2)}},$$

where $S(s) = 4s^3 - g_2 s - g_3$.

In her notation, $R_0(x_1, x_2) = -x_1^2 x_2^2 + 6\ell_1 x_1 x_2 + 2\ell c_0(x_1 + x_2) + c_0^2 - k^2$, and $P(x) = -x^4 + 6\ell_1 x^2 + 4\ell c_0 x + c_0^2 - k^2$. Then,

$$\begin{aligned}
R_0(x_1, x_2) &= -x_1^2 x_2^2 + 3\ell_1(x_1^2 + x_2^2) + 2\ell c_0(x_1 + x_2) + c_0^2 - k^2 - 3\ell_1(x_1 - x_2)^2 \\
&= R(x_1, x_2) - 3\ell_1(x_1 - x_2)^2
\end{aligned}$$

and so Kowalewski's constants are related to the constants of this paper as follows:

$$3\ell_1 = H, \quad c_0^2 - k^2 = \widetilde{K}_2, \quad 2\ell c_0 = -K_3 \ .$$

The change of variables, which has often been called "mysterious" in the subsequent literature, is nothing more than the mapping from $\mathcal{C} \times \mathcal{C}$ into $\Gamma \times \Gamma$ that assigns to $M = (x_1, u_1)$ and $N = (x_2, u_2)$ the values $P_1 = (s_1, S(s_1)) = M - N$ and $P_2 = (s_2, S(s_2)) = M + N$. The variables $s_1$ and $s_2$ in her paper coincide with $\xi_1$ and $\xi_2$ in this paper by the following calculation:

$$\begin{aligned}
s_1 &= \frac{R(x_1, x_2) - 3\ell_1(x_1^2 - x_2^2) - \sqrt{P(x_1)P(x_2)}}{2(x_1 - x_2)^2} + \frac{1}{2}\ell_1 \\
&= \frac{R(x_1, x_2) - \sqrt{P(x)P(y)}}{2(x_1 - x_2)^2} - \frac{3\ell_1}{2} + \frac{1}{2}\ell_1 \\
&= \frac{\theta_1}{2} - \ell_1 = \frac{2(\xi_1 + C)}{2} - \ell_1 = \xi_1
\end{aligned}$$

because $\ell_1 = C$. An analogous calculation shows that $s_2 = \xi_2$. Therefore, Kowalewski's change of variables is fundamentally a paraphrase of Theorem 7 in this paper.

The extremal equations in terms of the new coordinates $s_1$ and $s_2$ are obtained exactly in the same manner as in the original paper of Kowalewski. For convenience to the reader we include these calculations.

For the sake of continuity with the rest of the paper we shall use $(x, u)$, and $(y, v)$ to denote the points of $\mathcal{C}$, rather than $(x_1, u_1)$ and $(x_2, u_2)$, and use $(\xi_1, \eta_1)$ and $(\xi_2, \eta_2)$ instead of $(s_1, S(s_1))$ and $(s_2, S(s_2))$. Upon identifying $y$ with $\bar{x}$, we get that

$$-\frac{1}{4}\left(\frac{dy}{dt}\right)^2 = P(y) + \bar{q}(t)(x(t) - y(t))^2.$$

Since $\frac{d\xi_i}{\eta_i} = (-1)^i \frac{dx}{u} + \frac{dy}{v}$ with $i = 1, 2$, it follows that along each extremal curve $x(t)$

$$\left(\frac{d\xi_i}{dt}\right)^2 \frac{1}{\eta_i^2} = \frac{1}{u^2}\left(\frac{dx}{dt}\right)^2 + \frac{1}{v_2}\left(\frac{dy}{dt}\right)^2 + 2(-1)^i \frac{1}{uv}\frac{dx}{dt}\frac{dy}{dt}.$$

The substitution of

$$-\frac{1}{4}\left(\frac{dx}{dt}\right)^2 = u^2 + q(t)(x-y)^2 ,$$

$$-\frac{1}{4}\left(\frac{dy}{dt}\right)^2 7 = v^2 + \bar{q}(t)(x-y)^2 , \text{ and}$$

$$\frac{dx}{dt}\frac{dy}{dt} = R_0(x(t), y(t)) \quad \text{(Theorem 2)}$$

into the preceding expression leads to

$$-\frac{4}{\eta_i^2}\left(\frac{d\xi_i}{dt}\right)^2 = 2 + (x-y)^2\left(\frac{q(t)}{u^2(t)} + \frac{\bar{q}(t)}{v^2(t)}\right) + 2(-1)^i\frac{R_0(x,y)}{uv}$$

$$= \frac{2u^2v^2 + (x-y)^2(v^2q + u^2\bar{q}) + 2(-1)^iR_0(x,y)uv}{u^2v^2}$$

$$= \frac{2u^2v^2 - (x-y)^2(R_1 + K_4^2(x-y)^2) + 2(-1)^iR_0(x,y)uv}{u^2v^2}$$

because $v^2q + u^2\bar{q} + R_1 + K_4^2(x-y)^2 = 0$ (Theorem 1).

With the aid of $R_1(x-y)^2 = -R_0^2 + u^2v^2$ the above expression simplifies to

$$-\frac{4}{\eta_i^2}\left(\frac{d\xi_i}{dt}\right)^2 = \frac{u^2v^2 + R_0^2 + 2(-1)^iR_0uv - (x-y)^4K_4^2}{u^2v^2}$$

$$= \frac{(uv + (-1)^iR_0)^2 - (x-y)^4K_4^2}{u^2v^2}$$

$$= \frac{4(x-y)^4}{u^2v^2}\left(\left(\frac{R_0 + (-1)^2uv}{2(x-y)^2}\right)^2 - \frac{K_4^2}{4}\right) .$$

Recall that $\frac{R_0 + (-1)^iuv}{2(x-y)^2} = \xi_i - \frac{H}{6} + \frac{k}{2}$. Upon introducing two new constants

$$k_1 = -\frac{H}{6} + \frac{k}{2} + \frac{K_4}{2} \quad \text{and} \quad k_2 = -\frac{H}{6} + \frac{k}{2} - \frac{K_4}{2}$$

into the above equality we get

$$-\frac{1}{\eta_i^2}\left(\frac{d\xi_i}{dt}\right)^2 = \frac{(x-y)^4}{u^2v^2}(\xi_i - k_1)(\xi_i - k_2), \quad i = 1, 2.$$

Finally, $\xi_2 - \xi_1 = \frac{uv}{(x-y)^2}$, and hence

$$\left(\frac{d\xi_i}{dt}\right)^2 = \frac{U(\xi_i)}{(\xi_1 - \xi_2)^2} \quad \text{where} \quad U(\xi_i) = -(4\xi_i^3 - g_2\xi_i - g_3)(\xi_i - k_1)(\xi_i - k_2).$$

Consequently,

$$\frac{d\xi_1}{\sqrt{U(\xi_1)}} + \rho \frac{d\xi_2}{\sqrt{U(\xi_2)}} = 0 \quad \text{with} \quad \rho = \pm 1.$$

The last equation coincides with the celebrated equation of Kowalewski when $\rho = 1$.

## 5. The ratio of coefficients and meromorphic solutions

Kowalewski begins her investigations of the heavy top with a brilliant observation that both the top of Euler and the top of Lagrange, being integrable by means of elliptic integrals, admit Laurent series solutions of complex time for their general solution. In the papers listed in the bibliography she classifies all the cases of the heavy top whose solutions are meromorphic solutions of complex time, concluding that meromorphic solutions require that two of the coefficients must be equal to each other, say $c_1 = c_2$, in which case

( 1°) $c_1 = c_2 = c_3$, and $a_1, a_2, a_3$ arbitrary,

( 2°) $c_1 = c_2$, $c_3$ arbitrary, but $a_1 = a_2 = 0$,

( 3°) $c_1 = c_2 = 2c_3$, and $a_3 = 0$

are the only situations admitting general meromorphic solutions.

In this section we shall pursue Kowalewski's approach in a more general setting and classify all the cases in which the Hamiltonian solutions of

$$H = \frac{1}{2}\left(\frac{H_1^2}{c_1} + \frac{H_2^2}{c_2} + \frac{H_3^2}{c_3}\right) + a_1 h_1 + a_2 h_2 + a_3 h_3$$

are meromorphic functions of time under the assumption that $c_1 = c_2$.

We shall also include the following limiting cases in the subsequent analysis:

$$H = \frac{1}{2}H_3^2 + a_1 h_1 + a_2 h_2 + a_3 h_3$$

obtained as the common value $c_1 = c_2$ tends to $\infty$ while $c_3 = 1$, and

$$H = \frac{1}{2}(H_1^2 + H_2^2) + a_1 h_1 + a_2 h_2 + a_3 h_3$$

obtained as $c_3$ tends to $\infty$ while $c_1 = c_2 = 1$.

These limiting Hamiltonians are important for different reasons. The first case is completely integrable but falls outside of Kowalewski's classification, while the second case, which also falls outside of Kowalewski's classification, arises naturally in the theory of curves.

In what follows, we shall assume that all the variables $H_1, H_2, H_3, h_1, h_2, h_3$ are complex (but the coefficients $c_1, c_2, c_3$ and $a_1, a_2, a_3$ are real) and search for the cases where the solutions of

$$\text{(24)} \qquad \frac{dK}{dt} \;=\; [\Omega, K] + [a, P], \quad \frac{dP}{dt} \;=\; [\Omega, P] + k[a, K]$$

are meromorphic functions of complex time $t$, for at least an open set of initial values in $\mathbb{C}^6$.

It will be convenient to use $c$ instead of $c_3$, and to use $m$ to denote the common ratio $\frac{c_3}{c_1} = \frac{c_3}{c_2}$. That is, $m = \frac{c}{c_1}$ and hence $c_1 = \frac{c}{m} = c_2$. Therefore, $\widehat{\Omega} = \begin{pmatrix} \frac{m}{c} H_1 \\ \frac{m}{c} H_2 \\ \frac{1}{c} H_3 \end{pmatrix}$. It follows that $J_0 K J_0 = \Omega$ with

$$J_0 = \frac{1}{\sqrt{c}} \begin{pmatrix} 1 & 0 & 0 \\ 0 & 1 & 0 \\ 0 & 0 & m \end{pmatrix}.$$

The limiting case $H \;=\; \frac{1}{2} H_3^2 + a_1 h_1 + a_2 h_2 + a_3 h_3$ occurs when $m = 0$ and $c = 1$, while the other limiting case $H \;=\; \frac{1}{2}(H_1^2 + H_2^2) + a_1 h_1 + a_2 h_2 + a_3 h_3$, occurs when $\frac{m}{c} = 1$ and $\lim c = \infty$,

We shall use $\langle A, B \rangle$ to denote the trace of $-\frac{1}{2} AB$ for any antisymmetric matrices $A, B$ with complex entries (not to be confused with the Hermitian product on $\mathbb{C}^3$). Then,

$$\langle A, B \rangle \;=\; \widehat{A} \cdot \widehat{B} \;=\; a_1 b_1 + a_2 b_2 + a_3 b_3 \;.$$

In terms of this notation

$$H \;=\; \frac{1}{2} \langle J_0 K J_0, K \rangle + \langle a, P \rangle.$$

It is easy to verify that the Casimir functions

$$G \;=\; \langle P, P \rangle + k \langle K, K \rangle \quad \text{and} \quad J \;=\; \langle P, K \rangle$$

remain constants of motion for the complex system (24).

We shall now seek conditions on the ratio $m$ and the coefficients $a_1, a_2, a_3$ so that solutions $K(t), P(t)$ of equation (24) are of the form

$$K(t) \;=\; t^{-n_1}(K_0 + K_1 t + \cdots + K_n t^n + \cdots) \,,$$
$$P(t) \;=\; t^{-n_2}(P_0 + P_1 t + \cdots + P_n t^n + \cdots)$$

for some positive integers $n_1$ and $n_2$, and matrices $K_n, P_n$ in $\mathrm{so}_3(\mathbb{C})$ for each $n = 0, 1, \ldots$, with neither $K_0$ nor $P_0$ equal to zero.

Evidently not every solution has a pole in $\mathbb{C}$. For instance, $K = 0$ and $P = a$ is a solution for any choice of the constants in question. We shall assume that the meromorphic solutions occur for at least some open set of initial conditions in $\mathbb{C}^6$, which in turn implies that the meromorphic solutions should be parametrized by six arbitrary complex constants.

LEMMA 1. $n_1 = 1$ and $n_2 = 2$, provided that $[J_0 K_0 J_0, P_0] \neq 0$.

*Proof.* The leading terms in $\frac{dP}{dt} = [\Omega, P] + k[a, K]$ are given by $-n_2 t^{-(n_2+1)} P_0 = t^{-(n_2+n_1)}[J_0 K_0 J_0, P_0]$ which implies that $n_2 + 1 = n_1 + n_2$. Therefore, $n_1 = 1$.

Then the leading terms in $\frac{dK}{dt} = [\Omega, K] + [a, P]$ are given by

$$-t^{-2} K_0 = t^{-2}[J_0 K_0 J_0, K_0] + t^{-n_2}[a, P_0].$$

If $[a, P_0] = 0$ then $-K_0 = [J_0 K_0 J_0, K_0]$. The latter relation can hold only for $K_0 = 0$. Thus $[a, P_0] \neq 0$ and therefore $n_2 = 2$. The proof is now finished.

In the original paper Kowalewski claims the results of Lemma 1 without any proof ("on s'assure facilement, en comparant les exposants des premiers termes dans les membres gauches et dans les membres droits des équations considérées que l'on doit avoir $n_1 = 1, m_1 = 2$," [9, p. 178]). This claim, first criticized by A. A. Markov, is still considered an open gap in the original approach of Kowalewski. The assumption $[J_0 K_0 J_0, P_0] \neq 0$ is very likely inessential (for the contrary would reduce the number of possible solutions), and hence the original investigations of Kowalewski are in all probability correct in spite of the gap. We shall subsequently assume that $n_1 = 1$ and $n_2 = 2$ without going into details caused by vanishing of the above Lie brackets.

Upon equating the coefficients that correspond to the same powers of $t$ in equations (24) we come to the following relations:

$$(25) \qquad -K_0 = [J_0 K_0 J_0, K_0] + [a, P_0], \qquad 2P_0 = [J_0 K_0 J_0, P_0]$$

and

$$(26a) \qquad (n-1)K_n = \sum_{i=0}^{n}[J_0 K_i J_0, \ K_{n-i}] + [a, P_n] \ ,$$

$$(26b) \qquad (n-2)P_n = \sum_{i=0}^{n}[J_0 K_i J_0, \ P_{n-i}] + k[a, K_{n-2}]$$

for $n \geq 1$.

The same procedure applied to the Hamiltonian $H$ and the Casimir functions $G$ and $J$ gives

$$(27) \qquad \frac{1}{2}\sum_{i=0}^{n}\langle J_0 K_i J_0, \ K_{n-i}\rangle + \langle a, P_n\rangle = \delta_{2n}H \ ,$$

$$(28) \qquad \sum_{i=0}^{n}\langle P_i, P_{n-i}\rangle + k\sum_{i=0}^{n-2}\langle K_i, K_{n-2-i}\rangle = \delta_{4n}G \ ,$$

(29)
$$\sum_{i=0}^{n} \langle P_i, K_{n-i} \rangle \;=\; \delta_{3n} J \;,$$

with $\delta_{ij}$ denoting the Dirac function equal to 1 only for $i = j$, and otherwise equal to zero. It follows that

$$\sum_{i=0}^{2} \langle J_0 K_i J_0, K_{2-i} \rangle + \langle a, P_2 \rangle \;=\; H$$

$$\sum_{i=0}^{4} \langle P_i, P_{4-i} \rangle + k \sum_{i=0}^{2} \langle K_i, K_{2-i} \rangle \;=\; G$$

$$\sum_{i=0}^{3} \langle P_i, K_{3-i} \rangle \;=\; J$$

and consequently the first four stages of equations (26) must generate three arbitrary constants $H, G$ and $J$.

Following Kowalewski's original paper we shall also assume, without any loss in generality, that $a_2 = 0$ and that $a_1^2 + a_3^2 \neq 0$. This situation can always be arranged by a suitable rotation of $a$. We shall also rule out the cases $a_1 = 0$ as it corresponds to the case of Lagrange, and also rule out $m = 1$ as it corresponds to the well-known case $c_1 = c_2 = c_3$.

We shall now use $p_n, q_n, r_n$ to denote the entries of $K_n$ while $f_n, g_n, h_n$ will denote the entries of $P_n$ in the recursive relations (25) and (26). We then have:

THEOREM 1. *There are finitely many solutions of equation (25) if and only if $2m - 1 \neq 0$. They are given by*

(a)
$$p_0 \;=\; 0, \quad q_0 \;=\; \frac{2\varepsilon ic}{m}, \quad r_0 = 0, \quad f_0 = i\varepsilon h_0, \quad g_0 = 0, \quad h_0 = \frac{2c}{m(a_3 + \varepsilon i a_1)}$$

*and*

(b)
$$p_0 \;=\; -i\varepsilon g_0, \quad q_0(2m-1) \;=\; \frac{2a_3 c}{a_1}, \quad r_0 = 2\varepsilon ic,$$

$$f_0 \;=\; \frac{2c}{a_1}, \quad g_0 = i\varepsilon f_0, \quad h_0 = 0$$

*with $\varepsilon^2 = 1$ in either* (a) *or* (b).

*When $2m - 1 = 0$, then $a_3 = 0$ and $q_0$ is an arbitrary complex number corresponding to each case $\varepsilon = \pm 1$.*

*Proof.* The equation $2P_0 = [P_0, J_0 K_0 J_0]$ is linear in $P_0$ and can be written as

$$\text{(i)} \qquad \begin{pmatrix} 2 & -\frac{1}{c}r_0 & -\frac{m}{c}q_0 \\ \frac{1}{c}r_0 & 2 & \frac{m}{c}p_0 \\ \frac{m}{c}q_0 & -\frac{m}{c}p_0 & 2 \end{pmatrix} \begin{pmatrix} f_0 \\ g_0 \\ h_0 \end{pmatrix} = 0 .$$

The determinant of the preceding matrix is $2(4 + \frac{m^2}{c^2}(p_0^2 + q_0^2) + \frac{1}{c^2}r_0^2)$, and must vanish to admit nonzero solutions. Hence,

$$\text{(ii)} \qquad \frac{m^2}{c^2}(p_0^2 + q_0^2) + \frac{1}{c^2}r_0 = -4 .$$

It also follows from $2P_0 = [P_0, J_0 K_0 J_0]$ that $\langle P_0, P_0 \rangle = 0$ and that $\langle P_0, J_0 K_0 J_0 \rangle = 0$. Therefore,

$$\text{(iii)} \qquad f_0^2 + g_0^2 + h_0^2 = 0, \text{ and } m(p_0 f_0 + q_0 g_0) + r_0 h_0 = 0 .$$

Casimir relation (29) implies that $\langle P_0, K_0 \rangle = p_0 f_0 + q_0 g_0 + h_0 r_0 = 0$, which together with (iii) implies that $r_0 h_0 = 0$. Consider first the case $r_0 = 0$.

The relation $K_0 = [K_0, J_0 K_0 J_0] + [P_0, a]$ gives

$$\text{(iv)} \quad p_0 = \frac{(m-1)}{c}r_0 q_0 - a_3 g_0, \quad q_0 = \frac{(m-1)}{c}r_0 p_0 + a_3 f_0 - a_1 h_0, \quad r_0 = a_1 g_0 .$$

Hence $g_0 = 0$, and consequently $f_0^2 + h_0^2 = 0$. Thus, $f_0 = \varepsilon i h_0$ with $\varepsilon = \pm 1$.

It follows from (27) that $\frac{1}{2}\frac{m}{c}(p_0^2 + q_0^2) + \frac{1}{2c}r_0^2 + a_1 f_0 + a_3 h_0 = 0$. Thus,

$$(a_1 \varepsilon i + a_3) h_0 = a_1 f_0 + a_3 h_0 = -\frac{1}{2}\frac{m}{c}(p_0^2 + q_0^2) = \frac{-m^2 c(p_0^2 + q_0^2)}{2mc^2} = \frac{4c}{2m} = \frac{2c}{m}$$

($m \neq 0$, as can be seen from (ii)). Therefore,

$$h_0 = \frac{2c}{m(a_3 + \varepsilon i a_3)} .$$

Equations (iv) now imply that $p_0 = 0$, and

$$\begin{aligned} q_0 &= a_3 f_0 - a_1 h_0 = a_3 \varepsilon h_0 - a_1 h_0 = h_0(-a_1 + \varepsilon i a_3) \\ &= \frac{2}{m(a_3 + \varepsilon i a_1)}(-a_1 + \varepsilon i a_3) = \frac{2\varepsilon i c}{m}. \end{aligned}$$

We have now verified the solutions given by (a). The solutions in (b) correspond to $h_0 = 0$, and are obtained as follows:

Linear system (i) implies that $\frac{m}{c}(f_0 q_0 - g_0 p_0) = 0$. Combined with (iii) this relation gives that either $m = 0$, or that $p_0^2 + q_0^2 = 0$. In both cases $\frac{1}{c^2}r_0^2 = -4$ as can be seen from (ii). Therefore, $r_0 = 2\varepsilon i c$ with $\varepsilon = \pm 1$.

Then $g_0 = \frac{r_0}{a_1} = \frac{2\varepsilon i c}{a_1}$ (equations (iv)), and $f_0 = -\frac{r_0^2}{2a_1 c} = \frac{2c}{a_1}$, because $\frac{1}{2c}r_0^2 + a_1 f_0 + a_3 h_0 = 0$. That is, $g_0 = \varepsilon i f_0$.

The relation $p_0 f_0 + q_0 g_0 = 0$ gives $f_0(p_0 + \varepsilon i g_0) = 0$. Thus $p_0 = -\varepsilon i q_0$. Then,

$$q_0 = -\frac{(m-1)}{c} r_0 p_0 + a_3 f_0 = -\frac{(m-1)}{c}(2\varepsilon i)c(-\varepsilon i q_0) + a_3 f_0,$$

or

$$q_0(1 + 2(m-1)) = a_3 f_0 = \frac{2a_3 c}{a_1}.$$

Hence solutions in (b) are also verified concluding the proof of the theorem.

It remains to resolve the recursive relations (26). They are linear in $K_n$ and $P_n$ and can be expressed as

(30)

$$(n-1)K_n + [K_n, J_0 K_0 J_0] + [K_0, J_0 K_n J_0] + [P_n, a] = \sum_{i=1}^{n-1}[J_0 K_i J_0, K_{n-1-i}],$$

$$(n-2)P_n + [P_n, J_0 K_0 J_0] + [P_0, J_0 K_n J_0] = \sum_{i=1}^{n-1}[J_0 K_i J_0, P_{n-1-i}] + k[a, K_{n-1}]$$

for $n \geq 1$ with $K_{-1} = 0$.

This linear system in six variables $p_n, q_n, r_n, f_n, g_n, h_n$ can be written more explicitly as follows:

$$(n-1)p_n + \frac{1}{c}(mq_n r_0 - q_0 r_n + mq_0 r_n - q_n r_0) - a_3 g_n = A_n,$$

$$(n-1)q_n + \frac{1}{c}(r_n p_0 - mp_n r_0 + r_0 p_n - mp_0 r_n) + a_3 f_n - a_1 h_n = B_n,$$

$$(n-1)r_n + a_1 g_n = C_n$$

and

$$(n-2)f_n + \frac{1}{c}(mq_n h_0 - g_0 r_n + mq_0 h_n - g_n r_0) = D_n,$$

$$(n-2)g_n + \frac{1}{c}(r_n f_0 - mp_n h_0 + r_0 f_n - mp_0 h_n) = E_n,$$

$$(n-2)h_n + \frac{m}{c}(p_n g_0 - f_0 q_n + p_0 g_n - f_n g_0) = F_n$$

with $A_n, B_n, C_n, D_n, E_n, F_n$ denoting the appropriate quantities on the right-hand side of equation (30).

1. *Meromorphic solutions for $r_0 = 0$.* It follows from Lemma 2 that $p_0 = r_0 = g_0 = 0$. Therefore the preceding equations cascade into the following independent subsystems:

(31)

$$(n-1)p_n + \frac{(m-1)}{c}q_0 r_n - a_3 g_n = A_n,$$

$$(n-1)r_n + a_1 g_n = C_n,$$

$$(n-2)g_n + \frac{1}{c}(r_n f_0 - mp_n h_0) = E_n$$

and

(32)
$$(n-1)q_n + a_3 f_n - a_1 h_n \ = \ B_n \ ,$$
$$(n-2)f_n + \frac{1}{c}(mq_0 h_n + mq_n h_0) \ = \ D_n \ ,$$
$$(n-2)h_n - \frac{m}{c}(f_0 q_n + q_0 f_n) \ = \ F_n \ .$$

The determinant $\Delta$ of the overall system is the product of determinants $\Delta_1$ and $\Delta_2$ with

$$\Delta_1 = \ \det \begin{pmatrix} n-1 & \frac{(m-1)}{c}q_0 & -a_3 \\[2mm] 0 & n-1 & a_1 \\[2mm] -\frac{m}{c}h_0 & \frac{1}{c}f_0 & n-2 \end{pmatrix} \ ,$$

$$\Delta_2 = \ \det \begin{pmatrix} n-1 & a_3 & -a_1 \\[2mm] \frac{m}{c}h_0 & n-2 & \frac{m}{c}q_0 \\[2mm] -\frac{m}{c}f_0 & -\frac{m}{c}q_0 & n-2 \end{pmatrix} \ .$$

Then,

$$\Delta_2 \ = (n-1)[(n-2)^2 + \frac{m^2}{c^2}q_0^2] - \frac{m}{c}h_0((n-2)a_3 - a_1 \frac{m}{c}q_0)$$
$$- \ \frac{m}{c}f_0(a_3 \frac{m}{c}q_0 + a_1(n-2))$$
$$= (n-1)((n-2)^2 + \frac{m^2}{c^2}q_0^2) - (n-2)\frac{m}{c}(a_3 h_0 + f_0 a_1)$$
$$+ \ \frac{m^2}{c^2}q_0(a_1 h_0 - a_3 f_0).$$

Recalling that $q_0 = \frac{2\varepsilon ic}{m}$, $f_0 = i\varepsilon h_0$ and $h_0 = \frac{2c}{m(a_3 + i\varepsilon a_1)}$ we get that

$$\Delta_2 \ = \ (n-1)((n-2)^2 - 4) - 2(n-2) + 4 \ = \ (n+1)(n-2)(n-4).$$

Then,

$$\Delta_1 \ = (n-1)((n-1)(n-2) - \frac{a_1}{c}f_0) - \frac{m}{c}h_0(\frac{(m-1)}{c}\frac{a_1}{c}q_0 + a_3(n-1))$$
$$= (n-1)((n-1)(n-2) - a_1\frac{f_0}{c} - \frac{m}{c}h_0 a_3) - \frac{m(m-1)}{c}a_1 h_0 q_0$$
$$= (n-1)((n-1)(n-2) - \frac{2}{m}) - \frac{h_0}{c}(m(m-1)q_0 a_1 + (n-1)(m-1)a_3)$$
$$= (n-1)((n-1)(n-2) - \frac{2}{m}) - \frac{h_0}{c}(m-1)(2\varepsilon i a_1 - (n-1)a_3) \ .$$

But

$$\frac{h_0}{c}(2\varepsilon i a_1 - (n-1)a_3) = \frac{2(2\varepsilon i a_1 - (n-1)a_3)}{m(a_3 + \varepsilon i a_1)}$$

$$= \frac{2}{m(a_1^2 + a_3^2)}(2a_1^2 - (n-1)a_3^2 + i a_1 a_3(1+n)).$$

For $\Delta_1$ to be equal to zero the imaginary part of $\Delta_1$ must be equal to zero, which occurs only when $a_3 = 0$. But then

$$\Delta_1 = (n-1)((n-1)(n-2) - \frac{2}{m}) - \frac{4(m-1)}{m}$$

$$= \frac{(n-1)((n-1)(n-2)m - 2) - 4(m-1)}{m}$$

$$= \frac{(n-3)((n^2 - n + 2)m - 2)}{m}.$$

Therefore, $\Delta = 0$ when $n = 2$, $n = 4$, and $n = 3$ provided that $a_3 = 0$. In such a case $\Delta$ also vanishes for another positive integer value $k$ given by $(k^2 - k + 2)m - 2 = 0$.

THEOREM 2. *All solutions of equations* (26) *are parametrized by at most four arbitrary constants, and hence do not provide for general solutions of equation* (24).

*Proof.* The maximum number of solutions of equations (26) occurs when $n = 3$ is a zero of $\Delta$, that is when $a_3 = 0$. Our theorem follows from the simple observation that the kernel of the cascaded linear system with $a_3 = 0$,

$$(n-1)p_n + \left(\frac{m-1}{c}\right)q_0 r_n = 0, \qquad (n-1)q_n - a_1 h_n = 0,$$

$$(n-1)r_n + a_1 g_n = 0, \qquad (n-2)f_n + \frac{1}{c}mq_0 h_n + \frac{m}{c}q_n h_0 = 0,$$

$$(n-2)g_n + \frac{1}{c}(r_n f_0 - m h_0 p_n) = 0, \quad (n-2)h_n - \frac{m}{c}f_0 q_n - \frac{m}{c}q_0 f_n = 0,$$

is at most one dimensional at each singular value of $\Delta$. Recall that

$$f_0 = \varepsilon i h_0, \quad h_0 = \frac{2c}{ma_1} \quad \text{and that} \quad q_0 = \frac{2\varepsilon i c}{m}.$$

The foregoing linear systems are described by the following matrices:

$$\begin{pmatrix} n-1 & 0 & -a_1 \\ \frac{2}{a_1} & n-2 & 2\varepsilon i \\ -\frac{2\varepsilon i}{ma_1} & -2\varepsilon i & n-2 \end{pmatrix} \quad \text{and} \quad \begin{pmatrix} n-1 & 2\varepsilon i m(m-1) & 0 \\ 0 & n-1 & a_1 \\ -\frac{2}{a_1} & \frac{2\varepsilon i}{a_1} & n-2 \end{pmatrix}.$$

The rank of each of these matrices is at least two, and therefore there can be at most four constants arising from $n = 2$, $n = 3$, $n = 4$ and $(k^2 - k + 2)m - 2 = 0$. Our proof is finished because the number of parameters required for general solutions of equation (24) is six.

*Remark.* It might seem plausible that equations (26) do not admit solutions at the singular stages. Remarkably, that is not the case. In fact, for $n = 2$

$$\widehat{K}_2 = \begin{pmatrix} 0 \\ q_2 \\ 0 \end{pmatrix}, \quad \widehat{P}_2 = \begin{pmatrix} f_2 \\ 0 \\ h_2 \end{pmatrix} \quad \text{with}$$

$f_2 = \frac{\varepsilon i}{a_1} q_2 + \frac{k}{m} a_1$, and $h_2 = \frac{q_2}{a_1}$, with $q_2$ an arbitrary complex number. Then, for $n = 3$

$$\widehat{K}_3 = \begin{pmatrix} \frac{\varepsilon i(m-1)}{2m} a_1 g_3 \\[2mm] \frac{a_1}{2} h_3 \\[2mm] -\frac{a_1}{2} g_3 \end{pmatrix}, \quad \widehat{P}_3 = \begin{pmatrix} -i\varepsilon h_3 \\ g_3 \\ h_3 \end{pmatrix}, \quad \text{with} \quad h_3 = -\frac{1}{4} k a_1 q_2$$

and $g_3$ an arbitrary complex number. Finally, for $n = 4$

$$\widehat{K}_4 = \begin{pmatrix} 0 \\ q_4 \\ 0 \end{pmatrix}, \quad \widehat{P}_4 = \begin{pmatrix} \frac{-2n}{2a_1} q_2^2 - 2\varepsilon i q_4 \\ 0 \\ \frac{3}{a_1} q_4 \end{pmatrix}$$

with $q_4$ arbitrary.

2. *Meromorphic solutions for $h_0 = 0$.* In this case it is convenient to take $c = 1$. Then the solutions of equations (25) are given by

$$p_0 = -\varepsilon i q_0, \quad q_0(2m-1) = \frac{2a_3}{a_1}, \quad r_0 = 2\varepsilon i, \quad f_0 = \frac{2}{a_1}, \quad g_0 = i\varepsilon f_0, \quad h_0 = 0$$

with $\varepsilon = \pm 1$. We shall first suppose that $2m - 1 \neq 0$ and show that the number of constants that parametrize the general solution is less than six.

Since $p_0 + i\varepsilon q_0 = 0$ and $f_0 + \varepsilon i g_0 = 0$ it is natural to consider the following change of variables:

$$u_n = p_n + \varepsilon i q_n, \ v_n = p_n - \varepsilon i q_n, \ w_n = f_n + \varepsilon i q_n, \ z_n = f_n - i\varepsilon q_n.$$

It follows from equations (37a) and (37b) that

$$(n-1)u_n - i\varepsilon(m-1)r_0 u_n - i\varepsilon(m-1)u_0 r_n - i\varepsilon a_3 w_n - a_1 i\varepsilon h_n$$
$$= u_n(n-1+2(m-1)) - i\varepsilon a_3 w_n - i\varepsilon a_1 h_n = A_n + i\varepsilon B_n \ ,$$
$$v_n((n-1)-2(m-1)) + i\varepsilon(m-1)v_0 r_n + i\varepsilon a_3 z_n + i\varepsilon a_1 h_n = A_n - i\varepsilon B_n \ ,$$
$$(n-1)r_n + a_1 g_n = (n-1)r_n - \frac{\varepsilon i a_1 w_n}{2} + \frac{\varepsilon i q z_n}{2} = C_n \ ,$$
$$(n-2-2)w_n = D_n + i\varepsilon E_n \ ,$$
$$(n-2+2)z_n - \varepsilon i z_0 r_n + \varepsilon i m v_0 h_n = D_n + i\varepsilon E_n \ ,$$
$$(n-2)h_n + m(p_n g_0 - f_0 q_n + p_0 g_n - f_n q_0)$$
$$= (n-2)h_n + m f_0 \varepsilon i u_n - q_0 w_n = F_n \ .$$

The matrix corresponding to this linear system is given below:

$$M = \begin{pmatrix} 2m+n-3 & 0 & 0 & -i\varepsilon a_3 & 0 & -i\varepsilon a_1 \\ 0 & n+1-2m & i\varepsilon(m-1)v_0 & 0 & i\varepsilon a_1 & i\varepsilon a_1 \\ 0 & 0 & n-1 & \frac{-i\varepsilon a_1}{2} & \frac{i\varepsilon a_1}{2} & 0 \\ 0 & 0 & 0 & n-4 & 0 & 0 \\ 0 & 0 & -\varepsilon i z_0 & 0 & n & \varepsilon i m v_0 \\ m f_0 \varepsilon i & 0 & 0 & -q_0 & 0 & n-2 \end{pmatrix} .$$

The determinant $\Delta$ of the preceding matrix is given by $\Delta = (n-4)(n+1-2m)\Delta_1$ with $\Delta_1$ equal to the determinant of

$$\begin{pmatrix} n-3+2m & 0 & 0 & -i\varepsilon a_1 \\ 0 & n-1 & \frac{i\varepsilon a_1}{2} & 0 \\ 0 & -\varepsilon i z_0 & n & \varepsilon i m v_0 \\ m f_0 \varepsilon i & 0 & 0 & n-2 \end{pmatrix} .$$

Then

$$\begin{aligned}
\Delta_1 = {} & (n-3+2m)(n-2)(n(n-1) - \frac{a_1 z_0}{2}) \\
& - m f_0 \varepsilon i (-i\varepsilon a_1 (n(n-1) - \frac{a_1 z_0}{2})) \\
= {} & (n(n-1) - \frac{a_1 z_0}{2})((n-3+2m)(n-2) - ma_1 f_0) \\
= {} & (n(n-1) - 2)((n-3+2m)(n-2) - 2m) \\
= {} & (n(n-1) - 2)((n-1)(n-2) + 2(m-1)(n-2) - 2m) \\
= {} & (n+1)(n-2)((n-1)(n-2) + 2(m-1)(n-2) - 2m) \\
= {} & (n+1)(n-2)(n-3)(n-2+2m).
\end{aligned}$$

Therefore,

$$\Delta = (n+1)(n-2)(n-3)(n-4)(n+1-2m)(n-2+2m).$$

THEOREM 3. $2m - 1 = 0$ *and* $a_3 = 0$ *is the only case admitting a general meromorphic solution of equation* (24).

*Proof.* When $m \neq 0$ and $2m - 1 \neq 0$ then there are four constants that parametrize the solutions corresponding to $n = 2$, $n = 3$, $n = 4$ and one of the factors $n - 2 + 2m = 0$ or $n + 1 - 2m = 0$ because the kernel of $M$ is one dimensional for each singular value of $\Delta$. When $m = 0$, the kernel of $M$ is two dimensional for $n = 2$ and therefore contributes two constants that together with the constants corresponding to $n = 3$ and $n = 4$ again account for four constants.

In the remaining case $2m - 1 = 0$, and the zero-th stage (Theorem 1) accounts for two arbitrary constants, provided that $a_3 = 0$, that together with the constants corresponding to $n = 1$, $n = 2$, $n = 3$, $n = 4$ produce six arbitrary constants, which is the exact number required for the general solution.

This analysis shows that

(i)    $a = 0$,

(ii)   $c_1 = c_2$,      $a_1 = a_2 = 0$,

(iii)  $c_1 = c_2 = c_3$,   $a$ arbitrary,

(iv)   $c_1 = c_2 = 2c_3$,   $a_3 = 0$,

are the only cases that admit meromorphic solutions thus extending the conclusions of Kowalewski to a wider class of systems.

In the literature on Hamiltonian systems this claim of Kowalewski is often confused with the classification of completely integrable systems, that is, with systems that admit six independent integrals of motion all in involution with each other. This issue was considered by R. Liouville in 1896 ([11]) where he claimed that there is an extra integral of motion whenever the ratio $\frac{2c_3}{c_1}$ is an even integer. The exact nature of this claim is not altogether clear as the arguments in Liouville's paper are difficult to follow. In addition, the homogeneity properties of the equations for the heavy top, upon which any of Liouville's arguments are based, do not hold in our setting when $k \neq 0$ and further limit the relevance of Liouville's paper. However, the following simple example shows that there are additional cases of completely integrable systems which are outside of meromorphic class and therefore cannot be integrated by means of Abelian integrals.

3. *A completely integrable system not of Kowalewski type.* This case corresponds to the limiting ratio $m = 0$ and $a_3 = 0$. Then $H = \frac{1}{2}H_3^2 + a_1 h_1 + a_2 h_2$, and

$$\frac{dH_1}{dt} = H_2 H_3 - h_3 a_2, \quad \frac{dH_2}{dt} = -H_1 H_3 + h_3 a_1, \quad \frac{dH_3}{dt} = h_2 a_1 - h_2 a_1 \ ,$$

$$\frac{dh_1}{dt} = h_2 H_3 - k H_3 a_2, \quad \frac{dh_2}{dt} = -h_1 H_3 + k H_3 a_1, \quad \frac{dh_3}{dt} = k(H_1 a_2 - H_2 a_1) \ .$$

Let $w = h_1 + i h_2$. Then, $\frac{dw}{dt} = -i H_3(w - k(a_1 + i a_2))$. Hence, $\frac{d}{dt}(w - ka) = -i H_3(w - ka)$, with $a = a_1 + i a_2$, and consequently $|w - ka|^2$ is a constant of motion. As in all other integrable cases, this case is also completely integrable in the sense that there are six integrals of motion all in involution with each other.

Let $\theta$ be an angle defined by $w - ka = R e^{i\theta}$. Along each solution curve, $\frac{d}{dt}(w - ka) = R i e^{i\theta} \frac{d\theta}{dt} = -i H_3(R e^{i\theta})$ and hence $\frac{d\theta}{dt} = -H_3$. Then

$$\frac{1}{2}\left(\frac{d\theta}{dt}\right)^2 = \frac{H_3^2}{2} = H - R e \bar{a} w = H - R e \bar{a}(ka + R e^{i\theta})$$

$$= H - k|a|^2 - (R a_1 \cos\theta - R a_2 \sin\theta) \ .$$

So even though $\theta$ can be expressed in terms of elliptic integrals, it follows from the preceding results that the remaining equations cannot be integrated in terms of rational functions of $\theta$.

The other limiting case $H = \frac{1}{2}(H_1^2 + H_2^2) + a_1 h_1 + a_2 h_2 + a_3 h_3$ may be symplectically transformed into $\overline{H} = \frac{1}{2}(H_1^2 + H_3^2) + h_1$, provided that $a_2 = a_3 = 0$ and $a_1 = 1$. $\overline{H}$ corresponds to the total elastic energy of a curve $\gamma$ in $M$ given by $\frac{1}{2}\int(\kappa^2(t) + \tau^2(t))dt$ with $\kappa(t)$ and $\tau(t)$ denoting the curvature and the torsion of $\gamma$. While it is not known whether this Hamiltonian system admits an extra integral of motion, it nevertheless follows from the foregoing analysis that its equations cannot be integrated on Abelian varieties in terms of meromorphic functions.

*Acknowledgments.* I am grateful to Jean-Marie Strelcyn for pointing out the existence of Russian literature related to the gap in Kowalewski's paper concerning the order of poles in the solutions. I am also grateful to Ivan Kupka for his thoughtful comments on an earlier version of this paper.

University of Toronto, Toronto, ON, Canada
*E-mail address*: jurdj@math.toronto.edu